\documentclass{article}
\usepackage{graphicx}
\usepackage{amssymb}
\usepackage{amsmath}
\usepackage{color}
\usepackage{amsfonts}
\usepackage{changebar}
  \def\ZR{{\mathbb R}}

  \def\ZN{{\mathbb N}}

\def\beq{\begin{equation}}
\def\eeq{\end{equation}}
\def\be{\begin{equation}}
\def\ee{\end{equation}}
\def\beqar{\begin{eqnarray}}
\def\eeqar{\end{eqnarray}}
\def\ber{\begin{eqnarray}}
\def\eer{\end{eqnarray}}
\def\berb{\begin{eqnarray*}}
\def\eerb{\end{eqnarray*}}

\def\Ker#1{\mathop{\rm Ker}\nolimits#1}

\def\const{{\rm const}}

\def\norm#1.#2.{\|#1\|_{#2}}
\def\Norm#1.#2.{\big\|#1\big\|_{#2}}
\def\NOrm#1.#2.{\bigg\|#1\bigg\|_{#2}}
\def\NORm#1.#2.{\Big\|#1\Big\|_{#2}}
\def\NORM#1.#2.{\Bigg\|#1\Bigg\|_{#2}}

\newcommand{\vphi}{\varphi}

\def\vec#1{{\mathchoice{\mbox{\boldmath$\displaystyle#1$}}
{\mbox{\boldmath$\textstyle#1$}}
{\mbox{\boldmath$\scriptstyle#1$}}
{\mbox{\boldmath$\scriptscriptstyle#1$}}}}

\def \0b{{\hbox{\boldmath $0$}}}

 \newcommand{\bbb}{{\hbox{\bf b}}}

\newcommand{\eb}{\vec{e}} \newcommand{\fb}{\vec{f}}

 \newcommand{\nb}{\vec{n}}
 \newcommand{\pb}{\vec{p}}
\newcommand{\qb}{\vec{q}} 
 \newcommand{\tb}{\vec{t}}
\newcommand{\ub}{\vec{u}} \newcommand{\vb}{\vec{v}}
\newcommand{\wb}{\vec{w}} \newcommand{\xb}{\vec{x}}
\newcommand{\yb}{\vec{y}} \newcommand{\zb}{\vec{z}}

\newcommand{\Abb}{{\bf A}} \newcommand{\Bbbb}{{\bf B}}
 
\newcommand{\Ebb}{{\bf E}} 
 \newcommand{\Hbb}{{\bf H}}
\newcommand{\Ibb}{{\bf I}} 
\newcommand{\Kbb}{{\bf K}} 
\newcommand{\Mbb}{{\bf M}} 
 \newcommand{\Pbb}{{\bf P}}
\newcommand{\Qbb}{{\bf Q}} 
 
 \newcommand{\Vbb}{{\bf V}}

 \newcommand{\Fb}{\vec{F}}

 \newcommand{\Pb}{\vec{P}}
\newcommand{\Qb}{\vec{Q}} 
 
\newcommand{\Ub}{\vec{U}} \newcommand{\Vb}{\vec{V}}
\newcommand{\Wb}{\vec{W}}

\def \alb{\vec{\alpha}} \def \betab{\vec{\beta}}
\def \gamb{\vec{\gamma}} 
\def \deltab{\vec{\delta}} 
 
\def \thb{\vec{\theta}}

 \def \xib{\vec{\xi}}

 \def \phib{\vec{\phi}}
 
\def \psib{\vec{\psi}} \def \omegab{\vec{\omega}}
\def \Gamb{\vec{\Gamma}} 
\def \Thetab{\vec{\Theta}} 
\def \Xib{\vec{\Xi}} 
\def \Sib{\vec{\Sigma}} 
\def \Phib{\vec{\Phi}} 
\def \Omegab{\vec{\Omega}}

\newcommand{\cE}{{\cal E}} 
 
\newcommand{\cI}{{\cal I}} 
 
 \newcommand{\cN}{{\cal N}}

 \newcommand{\cV}{{\cal V}}

\def \ptb{\vec{\tilde{p}}}
\def \qtb{\vec{\tilde{q}}}

\newcounter{primjer}[section]
\newcounter{tvrdnja}[section]
\newcounter{zadatak}[section]
\newcommand {\proof} {\par{\it Proof}. \ignorespaces}
\newcommand {\eproof}

\newcommand {\mat}      [1] {\left[\begin{array}{#1}}
\newcommand {\rix}          {\end{array}\right]}
\long\def\Hidden#1{\relax}

\newtheorem{theorem}{Theorem}
\newtheorem{remark}{Remark}
\newtheorem{lemma}{Lemma}

\newtheorem{corollary}{Corollary}

\newcommand{\bbP}{{\mathbb P}}

\newcommand{\Pbbb}{{\mathbb P}}

\renewcommand{\Im}{\mathop{\rm Im}}

\newcommand{\thetab}{\boldsymbol{\theta}}

\newcommand{\diag}{\rm {diag}}

\def \Omegab{\vec{\Omega}}

\newcommand{\funv}{\psib}
\newcommand{\funu}{\phib}
\newcommand{\funS}{\Sib}
\newcommand{\funG}{\Gamb}
\title{Modeling and discretization methods for the
numerical simulation of elastic stents
\thanks{%
Received... Accepted... Published online on... Recommended by....
This work has been supported by {\em Deutscher Akademischer Austauschdienst} (DAAD)  via Project {\em Asymptotic and algebraic analysis of nonlinear eigenvalue problems in contact mechanics and electro
magnetism}. The third author is also supported by {\em Einstein Foundation Berlin via Einstein Center ECMath Project}: {\em Model Reduction for Nonlinear Parameter-Dependent Eigenvalue Problems in Photonic Crystals}.}}

\author{Luka Grubi\v{s}i\'{c}\footnotemark[2]
        \and Matko Ljulj\footnotemark[2]
        \and Volker Mehrmann\footnotemark[3]
        \and Josip Tamba\v{c}a\footnotemark[2]}

\begin{document}

\maketitle

\renewcommand{\thefootnote}{\fnsymbol{footnote}}

\footnotetext[2]{Department of Mathematics, Faculty of Science, University of Zagreb, Bijeni\v{c}ka 30, 10000 Zagreb, \texttt{\{luka,mljulj,tambaca\}@math.hr}.}
%\footnotetext[3]{Department of Mathematics, Faculty of Science, University of Zagreb, Bijeni\v{c}ka 30, 10000 Zagreb, \texttt{mljulj@math.hr}.}
\footnotetext[4]{Institut f\"ur Mathematik MA 4-5, TU Berlin, Str. des 17. Juni 136,
D-10623 Berlin, FRG. \texttt{mehrmann@math.tu-berlin.de}.}
%\footnotetext[5]{Department of Mathematics, Faculty of Science, University of Zagreb, Bijeni\v{c}ka 30, 10000 Zagreb, \texttt{tambaca@math.hr}.}

\begin{abstract}
A new model description for the numerical simulation of elastic stents is proposed. Based on the new formulation an $\inf$-$\sup$ inequality for the finite element discretization is proved and the proof of the $\inf$-$\sup$ inequality for the continuous problem is simplified. The new formulation also leads to faster simulation times despite an increased number of variables. The techniques also simplify the analysis and numerical solution of the evolution problem describing the movement of the stent under external forces. The results are illustrated via numerical examples.
\end{abstract}


\textbf{Ams Subj. Classification:}
74S05, 74K10, 74K30, 74G15, 74H15, 	65M15, 65M60

%\tableofcontents
\section{Introduction}
In this paper we present a new model description for the dynamic and stationary simulation of stents. This new formulation is using constrained partial differential equations in mixed variational weak form, which are based on a network structure consisting of one dimensional curved rods (struts). The new formulation will turn out to be particularly
convenient for the analysis of the partial differential equation, in particular in proving an $\inf$-$\sup$ inequality which directly transfers to a discrete $\inf$-$\sup$ inequality in the discretized setting, so that from classical results of \cite{BBF} the error estimates follow.

In the new formulation, the inextensibility and unshearability
of the rod are expressed in the weak formulation, but the continuity of the displacement and the modeling of infinitesimal rotations are modelled via constraints so that they do not have to be incorporated in  the function spaces as was done in the classical approach in \cite{RadHAZU}. This advantage of the new formulation comes along with the introduction of new unknowns for the displacements and infinitesimal rotations at vertices where different struts are connected and further unknowns for the contact couples and contact forces at the end points of each strut.
However, despite the introduction of many new unknowns, the numerical solvers become more efficient for the large scale cases.

%Furthermore, when we formulate the discrete projection of the mixed formulation, it turns %In order to do that we will
%express all relevant quantities directly in the mixed formulation as opposed to putting %them in the
%structure of the space of functions on which we formulate the variational forms as was %done in \cite{RadHAZU}.
%\subsection{Geometry of the stent and notation}\label{sec:geometry}
%\marginpar{do not have a Palmaz stent picture}
\begin{figure}[h!]
\begin{center}
\includegraphics[bb = 0 0 452 106, width=0.5\textwidth]{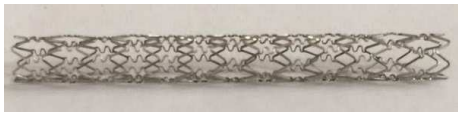}
\caption{Example of an elastic stent (Cypher stent by Cordis Corporation)}\label{cypher_small}
\end{center}
\end{figure}
%\bar{}
Stents, see Figure~\ref{cypher_small},
are typically considered as a union of struts each of which is modeled by a 1D curved rod model, see \cite{JT1,JT2}, and a set of junction conditions describing the connection of the struts,  see \cite{Griso}. The so-obtained model describes the three-dimensional behavior of stents but it has the complexity of a one-dimensional model. This model can be applied for any elastic structure made of thin curved (or straight) rods. It was first formulated in \cite{SIAMstent} and then reformulated in the weak form in \cite{IMASJ}. The properties of the mixed formulation for the model have  been analyzed in \cite{RadHAZU}, numerical methods have been introduced, and error estimates have been derived in \cite{GT}, but up to now error estimates for the contact forces, which are represented as Lagrange multipliers of the contact conditions, were missing. To derive these error estimates is one of the main results of this paper.

To model the topology of the stent we use an undirected graph  $\cN = (\cV, \cE)$ consisting of a set
$\cV$ of $n_\cV$ vertices, which are the points where the middle lines of the rods meet and a set $\cE$ of $n_\cE$ edges that represent a 1D description of the curved rod. To be able to use a 1D curved rod model, we additionally need to prescribe the local geometry of the rod, i.e., the middle curve and the geometry of the cross-section as well as the material properties of the stent. These are given by
\begin{itemize}
\item the function $\Phib^i:[0,\ell^i] \to \ZR^3$ as natural parametrization of the middle line of the  $i$th strut of length $\ell^i$, represented by the edge $\eb^i \in \cE$,
\item the shear modulus $\mu_i$ and the Young modulus $E_i$ as parameters describing the material of the $i$th strut,
\item as well as the width $w^i$ and the thickness $t^i$ of the rectangular cross-section of the $i$th strut.
\end{itemize}
Using  these quantities, in the stationary case, see \cite{SIAMstent}, the model for the $i$th strut $e_i\in \cE$, is given by the following system of ordinary differential equations (in space)
\begin{eqnarray}
\label{lj1ri} &&0 = \partial_s \pb^i + \fb^i,\\
\label{lj2ri} &&0 = \partial_s \qb^i  +  \tb^i \times\pb^i%
%+ \lb
,\\
\label{lj3ri} &&0 = \partial_s \omegab^i - \Qbb^i (\Hbb^i)^{-1} (\Qbb^i)^T \qb^i,\\
\label{lj4ri}
&&0 = \partial_s \ub^i +\tb^i \times  \omegab^i.
\end{eqnarray}
where for the $i$th strut
\begin{itemize}
\item $\ub^i:[0,\ell^i] \to \ZR^3$ denotes the vector of displacements on the middle curve,
\item     $\omegab^i:[0,\ell^i] \to \ZR^3$ is the vector of infinitesimal rotations of the cross-section,
\item $\qb^i$ is the contact moment and $\pb^i$ is the contact force,
\item  $\fb^i$ is the line density of the applied forces,
\item $\Qbb^i= [\tb^i, \nb^i, \bbb^i]$ is an orthogonal  rotation matrix associated to the middle curve, with $\tb^i=(\Phib^i)'$ being the unit tangent to the middle curve and $\nb^i$, $\bbb^i$ being vectors spanning the normal plane to the middle curve, so that $\Qbb^i$ represents the local basis at each point of the middle curve,
\item $\Hbb^i=\diag(\mu^i K^i, E^i I_n^i, E^i I_b^i)$ is a positive definite diagonal matrix, with the Young modulus $E^i$, the shear modulus $\mu^i$, $I_n^i$, $I_b^i$ are the moments of inertia of the cross section and $\mu^iK^i$ is the torsional rigidity of the cross section.
\end{itemize}
Equations (\ref{lj1ri})  and (\ref{lj2ri}) represent equilibrium equations (for forces and moments), while (\ref{lj3ri}) and (\ref{lj4ri}) are constitutive relations. In particular, (\ref{lj4ri}) describes the inextensibility and unshearability of the struts, see \cite{IMASJ} for more details.

In addition to equations (\ref{lj1ri})--(\ref{lj4ri}), at each vertex of the stent we have a kinematic coupling condition that $\ub$ and $\omegab$ are continuous and a  dynamic coupling condition describing the  balance of contact forces $\pb$ and contact moments $\qb$.

Denoting by $J^-_j$  the set of all edges that leave the $j$th vertex, i.e., the local variable is  equal to $0$ at vertex $j$ and by $J^+_j$  the set of all edges that enter the vertex, i.e.,  the local variable is  equal to $\ell^i$ for $i$th edge at the vertex $j$.
With these notations we obtain the node conditions
\begin{equation}\label{ccs}
\aligned
%\Omegab^j =
\omegab^i(0) &= \omegab^k(\ell^k), \qquad i \in J^-_j, k\in J^+_j, \qquad &j=1,\ldots, n_\cV,\\
%\Ub^j =
\ub^i(0) &= \ub^k(\ell^k), \qquad i \in J^-_j, k\in J^+_j, \qquad &j=1,\ldots, n_\cV,\\
\sum_{i\in J^+_j} \pb^i(\ell^i) - \sum_{i\in J^-_j} \pb^i(0) &= 0, \qquad &j=1,\ldots, n_\cV,\\
\sum_{i\in J^+_j} \qb^i(\ell^i) - \sum_{i\in J^-_j} \qb^i(0) &= 0, \qquad &j=1,\ldots, n_\cV.
\endaligned
\end{equation}
Since this is a pure traction problem, we can integrate over $s \in [0,\ell^i]$ and  specify a unique solution by requiring the two additional conditions
\begin{equation}\label{fix}
\int_\cN \ub := \sum_{i=1}^{n_\cV}\int_0^{\ell^i} \ub^i \, ds = 0, \qquad \int_\cN \omegab := \sum_{i=1}^{n_\cV}\int_0^{\ell^i} \omegab^i \, ds=0,
\end{equation}
which means that the total displacement as well as the total infinitesimal rotation are zero.

In the formulation of \cite{RadHAZU}, the model is described on the collection of all %$H^1$
displacements $\ub^i$ and infinitesimal rotations $\omegab^i$ for all edges which are continuous on the whole stent. Thus, the tuples of unknowns in the problem $\ub_S=((\ub^1,\omegab^1),\ldots,(\ub^{n_\cE},\omegab^{n_\cE}))$ belong to the space %$H^1(\cN;\ZR^6)$, where
\[
\aligned
H^1(\cN; \ZR^k) = \bigg\{& (\yb^1, \ldots, \yb^{n_\cE}) \in \prod_{i=1}^{n_\cE} H^1(0,\ell^i; \ZR^k) :\\
& \yb^i(0) = \yb^k(\ell^k), i \in J^-_j, k \in J^+_j, j=1,\ldots,n_\cV \bigg\},
\endaligned
\]
with $H^1(0,\ell^i; \ZR^k)$ being the Sobolov space of functions on $[0,\ell^i]$ whose derivatives up to the first derivative are square Lebesgue integrable.
%We also set the space of multipliers
%$$
%Q_S = L^2(\cN;\ZR^3) \times \ZR^3 \times \ZR^3= \prod_{i=1}^{n_\cE} L^2 (0,\ell^i; \ZR^3) \times \ZR^3 \times \ZR^3
%$$
The formulation given in \cite{RadHAZU} is a mixed formulation with Lagrange multipliers appearing in the formulation due to the inextensibility and unshearability of the struts in the 1d curved rod model (\ref{lj4ri}), and the two conditions on the total displacement and infinitesimal rotation (\ref{fix}). Under these conditions, in \cite{RadHAZU} an $\inf$-$\sup$ inequality was proved and the well-posedness of the problem was established. However, for the discretized  problem via the finite element method it would be necessary to also have a discrete $\inf$-$\sup$ inequality to obtain an error estimate also for the Lagrange multipliers approximation.

We will show that the elastic energy is coercive in the space implementing inextensibility, unshearability, and continuity of displacement and infinitesimal rotation. Using the results of
\cite{BBF}, we will show that the Lagrange multipliers are unique in both the old and the new formulation and for the discrete problem also in the new formulation.
This is a considerable improvement compared to the classical mixed formulation \cite{GT}, where a proof of the discrete $\inf$-$\sup$ inequality was not successful.

Let $\Abb_\cI \in \ZR^{3n_\cV, 3n_\cE}$ denote the incidence matrix of the oriented graph $(\cV,\cE)$ with three connected components, organized in the following way: a $3\times3$ submatrix at rows $3i-2,3i-1,3i$ and columns $3j-2,3j-1,3j$ is $\Ibb_3$ if the edge $j$ enters the vertex $i$, $-\Ibb_3$ if it leaves the vertex $i$ or $0$ otherwise. Then the matrix
$\Abb^+_\cI$ is obtained from $\Abb_\cI$ by setting all elements  $-1$ to $0$ and  $\Abb^-_\cI$ is defined as $\Abb_\cI=\Abb^+_\cI-\Abb^-_\cI$.
Let us also introduce the projectors
\[
\Pbbb^i_\cE \in \ZR^{3,3n_\cE}, \qquad \Pbbb^j_\cV\in \ZR^{3,3n_\cV}
\]
on the coordinates $3i-2,3i-1,3i$ and $3j-2,3j-1,3j$, respectively.
We will also need the spaces
\[
L^2(\cN;\ZR^3) = \bigotimes_{i=1}^{n_\cE} L^2(0,\ell^i; \ZR^3), \qquad L^2_{H^r}(\cN;\ZR^3) = \bigotimes_{i=1}^{n_\cE} {H^r}(0,\ell^i; \ZR^3),\quad {r\geq 1}
\]
with associated norms
\[
\aligned
&\|(\yb^1,\ldots,\yb^{n_\cE})\|_{L^2(\cN;\ZR^3)} = \left( \sum_{i=1}^{n_\cE} \|\yb^i\|_{L^2(0,\ell^i;\ZR^3)}^2\right)^{1/2},\\
&\|(\yb^1,\ldots,\yb^{n_\cE})\|_{L^2_{H^r}(\cN;\ZR^3)} = \left( \sum_{i=1}^{n_\cE} \|\yb^i\|_{H^{r}(0,\ell^i;\ZR^3)}^2\right)^{1/2}.
\endaligned
\]
The norm corresponding to the last term for $r=1$ is also used as the norm for $H^1(\cN;\ZR^6)$.
For a function $\yb=(\yb^1, \ldots,\yb^{n_\cE})\in L^2_{H^1}(\cN;\ZR^3)$ by $\yb'$ we denote $(\partial_s \yb^1, \ldots, \partial_s \yb^{n_\cE}) \in L^2(\cN;\ZR^3)$.

The results that we prove for the continuous model in Section~\ref{S4} hold for general geometries, while the results for the discrete approximation in Section~\ref{S5} are proved only for stent geometries with straight struts.

The paper is organized as follows. In Section~\ref{newmodel} we present the new formulation of the stent model. In Section~\ref{S4} we analyze the new model and give a proof of the inf-sup inequality for the weak formulation of the continuous infinite dimensional model and in Section~\ref{S5} we analyze the discrete model that is obtained after finite element discretization and show a corresponding inf-sup inequality. Finally in Section~\ref{evol} we study the dynamical system of the stent movement under excitation forces. We analyze the properties and present numerical simulation results.

\section{New formulation of the model}\label{newmodel}
\setcounter{equation}{0}

In this section we reformulate the stent model in a way that enables the proof of  a discrete $\inf$-$\sup$ inequality and then, using classical results, appropriate error estimates follow. For this, we treat all unknowns in the problem explicitly and do not encode them in the function spaces or the weak formulation.
In this way not only the inextensibility and unshearability of the rod is expressed in the weak formulation, but the continuity of the displacement and infinitesimal rotation is reflected in the function space of the mixed formulation in $H^1(\cN;\ZR^6)$. This leads to the introduction of new unknowns, displacements and infinitesimal rotations at vertices, and further, the contact moments and contact forces at the ends of each strut.

%One-dimensional equilibrium model for curved elastic rods we use here is given by the following first order system. For a given force with line density $\fb$ the model reads: find $(\ub,\omegab,\qb,\pb)$ such that
%\begin{align}
%\label{lj1r} 0 &= \partial_s \pb + \fb,\\
%\label{lj2r} 0 &= \partial_s \qb  +  \tb \times\pb%
%%+ \lb
%,\\
%\label{lj3r} 0 &= \partial_s \omegab - \Qbb \Hbb^{-1} \Qbb^T \qb,\\
%\label{lj4r}
%0 &= \partial_s \ub +\tb \times  \omegab.
%\end{align}

Since $\ub$ and $\omegab$ are continuous over the whole stent, we introduce as extra variables the displacements and infinitesimal rotations at the vertices $\Ub^i, \Omegab^i$,  $i=1,\ldots,n_\cV$, and then form the vectors
\[
\Ub=[\Ub^1, \ldots,\Ub^{n_\cV}]^T, \ \Omegab=[\Omegab^1, \ldots,\Omegab^{n_\cV}]^T.
\]
%
%We define $\Ub=(\Ub^1, \ldots,\Ub^{n_\cE})$ and $\Omegab=(\Omegab^1, \ldots,\Omegab^{n_\cE})$
Then the kinematic coupling at the vertex $j$ leads to the conditions
\begin{equation}\label{displ}
\aligned
&\ub^i(\ell^i) = \Ub^j, \ i \in J^+_j, \quad \ub^i(0) = \Ub^j, \ i \in J^-_j,\\
&\omegab^i(\ell^i) = \Omegab^j, \  i \in J^+_j, \quad  \omegab^i(0) = \Omegab^j, \  i \in J^-_j.
\endaligned
\end{equation}
To express the dynamic coupling conditions, we introduce the contact moments and forces at the ends of the struts,
\begin{equation}\label{contact_forces}
\Qb^i_+ = \qb^i(\ell^i), \quad \Qb^i_- = \qb^i(0), \quad \Pb^i_+ = \pb^i(\ell^i), \quad \Pb^i_- = \pb^i(0), \qquad i=1,\ldots,n_\cE
\end{equation}
and define
\[
\Pb_\pm = (\Pb^1_\pm, \ldots, \Pb^{n_\cE}_\pm), \qquad \Qb_\pm = (\Qb^1_\pm, \ldots, \Qb^{n_\cE}_\pm).
\]
Then the dynamic coupling conditions at the vertex $j$ can be expressed as
\begin{equation}\label{sila}
\sum_{i\in J^+_j} \Pb^i_+ -\sum_{i\in J^-_j} \Pb^i_- =0, \quad\sum_{i\in J^+_j} \Qb^i_+ -\sum_{i\in J^-_j} \Qb^i_- =0, \qquad j=1,\ldots, n_\cE.
\end{equation}
%
%Also, to fix the unique solution we add the conditions
%\begin{equation}\label{fix}
%\int_\cN \ub = \int_\cN \omegab =0.
%\end{equation}
Equations (\ref{lj1ri})--(\ref{lj4ri}), (\ref{displ}), (\ref{contact_forces}), (\ref{sila}), and (\ref{fix}) together constitute the stent problem in our new formulation for which we now derive in detail the weak formulation.
%Since the formulation is new we do all steps of the derivation in detail.

We multiply the $i$th equation of (\ref{lj1ri}) by $\vb^i\in H^1(0,\ell^i;\ZR^3)$ and that of (\ref{lj2ri}) by $\wb^i\in H^1(0,\ell^i;\ZR^3)$, add them, integrate over $s\in [0,\ell^i]$, and sum the equations over $i$. This yields
\[
\aligned
0 = \sum_{i=1}^{n_\cE} \int_0^{\ell^i} \partial_s \pb^i \cdot \vb^i + \fb^i \cdot \vb^i + \partial_s \qb^i  \cdot \wb^i+  \tb^i \times\pb^i \cdot \wb^i\, ds.
\endaligned
\]
After partial integration we obtain
\[
\aligned
0 = \sum_{i=1}^{n_\cE} \int_0^{\ell^i} \left(- \pb^i \cdot \partial_s \vb^i + \fb^i \cdot \vb^i -  \qb^i  \cdot \partial_s \wb^i {-}  \tb^i \times \wb^i \cdot\pb^i \right)\, ds + \pb^i \cdot \vb^i \big|_0^{\ell^i}+ \qb^i  \cdot \wb^i \big|_{0}^{\ell^i},
\endaligned
\]
i.e.,
\begin{equation}\label{lj12}
\aligned
\sum_{i=1}^{n_\cE} &\int_0^{\ell^i} \left(- \pb^i \cdot (\partial_s \vb^i +\tb^i \times \wb^i)  -  \qb^i  \cdot \partial_s \wb^i \right)\, ds\\
 &+ \Pb^i_+ \cdot \vb^i(\ell^i) - \Pb^i_- \cdot \vb^i(0)+ \Qb^i_+ \cdot \wb^i(\ell^i) - \Qb^i_- \cdot \wb^i(0)= - \sum_{i=1}^{n_\cE} \int_0^{\ell^i} \fb^i \cdot \vb^i\, ds.
\endaligned
\end{equation}
In a similar way we multiply (\ref{lj3ri}) by $\xib^i \in L^2(0,\ell^i;\ZR^3)$ and (\ref{lj4ri}) by $\thetab^i \in L^2(0,\ell^i;\ZR^3)$ integrate over $s\in [0,\ell^i]$ and sum all equations to obtain
\begin{equation}\label{lj34}
0 = \sum_{i=1}^{n_\cE}\int_0^{\ell^i} -\partial_s \omegab^i \cdot \xib^i + \Qbb^i (\Hbb^i)^{-1} (\Qbb^i)^T \qb^i \cdot \xib^i -  (\partial_s \ub^i +\tb^i \times  \omegab^i) \cdot \thetab^i\, ds.
\end{equation}
We also multiply the equations in (\ref{sila}) for the $j$th vertex by $\Vb^j$ and $\Wb^j$ from $\ZR^3$, respectively, and sum the equations over $j$ which gives
\[
\sum_{j=1}^{n_\cV} \left(\sum_{i\in J^+_j} \Pb^i_+ -\sum_{i\in J^-_j} \Pb^i_- \right)\cdot \Vb^j + \sum_{j=1}^{n_\cV} \left(\sum_{i\in J^+_j}\Qb^i_+ -\sum_{i\in J^-_j} \Qb^i_- \right)\cdot \Wb^j=0.
\]
Since $\sum_{i\in J^+_j} \Pb^i_+ = \mathbb{P}^j_\cV \Abb^+_\cI \Pb_+$ and $\sum_{i\in J^-_j} \Pb^i_- = \mathbb{P}^j_\cV \Abb^-_\cI \Pb_-$, this equation can be written as
\[
\sum_{j=1}^{n_\cV} \mathbb{P}^j_\cV \left(\Abb^+_\cI \Pb_+ -\Abb^-_\cI \Pb_- \right)\cdot \Vb^j +
\sum_{j=1}^{n_\cV} \mathbb{P}^j_\cV \left(\Abb^+_\cI \Qb_+ - \Abb^-_\cI \Qb_-\right)\cdot \Wb^j=0,
\]
and, therefore,
\begin{equation}\label{sila2}
\left(\Abb^+_\cI \Pb_+ -\Abb^-_\cI \Pb_- \right)\cdot \Vb + \left(\Abb^+_\cI \Qb_+ - \Abb^-_\cI \Qb_-\right)\cdot \Wb=0 %\qquad \Vb,\Wb \in \ZR^3.
\end{equation}
for all $\Vb=[\Vb^1, \ldots, \Vb^{n_\cV}]^T$, $\Wb=[\Wb^1, \ldots, \Wb^{n_\cV}]^T \in \ZR^{3n_\cV}$.

Multiplying the equations for the displacements in (\ref{displ}) by $\Thetab^i_+$ and $\Thetab^i_-$, we obtain
\[
\ub^i(\ell^i) \cdot \Thetab^i_+ = \Ub^j \cdot \Thetab^i_+, \quad i \in J^+_j, \qquad \ub^i(0) \cdot \Thetab^i_-= \Ub^j \cdot \Thetab^i_-, \quad i \in J^-_j.
\]
Since  $\mathbb{P}^i_\cE (\Abb^+_\cI)^T\Ub=\Ub^j$ for  $i \in J^+_j$, for $\Thetab_\pm=[\Thetab^1_\pm, \ldots, \Thetab^{n_\cE}_\pm]^T$ we have
%$$
%\Ub^j \cdot \Thetab^i_+ = \mathbb{P}^i_\cE (\Abb^+_\cI)^T\Ub \cdot \Thetab^i_+, \qquad \mbox{ for } i \on J^+_j.
%$$
%Thus we have
\[
\sum_{i=1}^{n_\cE} \ub^i(\ell^i) \cdot \Thetab^i_+ = (\Abb^+_\cI)^T\Ub \cdot \Thetab_+, \quad \sum_{i=1}^{n_\cE} \ub^i(0) \cdot \Thetab^i_- = (\Abb^-_\cI)^T\Ub \cdot \Thetab_-, \qquad \Thetab_+, \Thetab_- \in \ZR^{3n_\cE},
\]
and similarly, for the rotations, using the notation $\Xib_\pm=[\Xib^1_\pm, \ldots, \Xib^{n_\cE}_\pm]^T$, we get
\[
\sum_{i=1}^{n_\cE} \omegab^i(\ell^i) \cdot \Xib^i_+ = (\Abb^+_\cI)^T\Omegab \cdot \Xib_+, \quad \sum_{i=1}^{n_\cE} \omegab^i(0) \cdot \Xib^i_- = (\Abb^-_\cI)^T\Omegab \cdot \Xib_-,  \qquad \Xib_+, \Xib_- \in \ZR^{3n_\cE}.
\]
Thus, we have
\begin{equation}\label{displ2}
\sum_{i=1}^{n_\cE} (\ub^i(\ell^i) \cdot \Thetab^i_+ - \ub^i(0) \cdot \Thetab^i_- ) - (\Abb^+_\cI)^T\Ub \cdot \Thetab_+ + (\Abb^-_\cI)^T\Ub \cdot \Thetab_-=0, \qquad \Thetab_+, \Thetab_- \in \ZR^{3n_\cE},
\end{equation}
for the displacements and
\begin{equation}\label{rot2}
\sum_{i=1}^{n_\cE} (\omegab^i(\ell^i) \cdot \Xib^i_+ - \omegab^i(0) \cdot \Xib^i_- ) - (\Abb^+_\cI)^T\Omegab \cdot \Xib_+ + (\Abb^-_\cI)^T\Omegab \cdot \Xib_-=0, \qquad \Xib_+, \Xib_- \in \ZR^{3n_\cE}
\end{equation}
for the rotations.
We multiply the equations (\ref{fix}) by $\alb$ and $\betab$, respectively, and summing up, we obtain
\begin{equation}\label{*}
\alb \cdot \int_\cN \ub + \betab \cdot \int_\cN \omegab =0, \qquad \alb, \betab \in \ZR^3.
\end{equation}
Subtracting (\ref{sila2}) from (\ref{lj12}), we obtain
\begin{equation}\label{lj12A}
\aligned
\sum_{i=1}^{n_\cE} &\int_0^{\ell^i} \left(- \pb^i \cdot (\partial_s \vb^i +\tb^i \times \wb^i)  -  \qb^i  \cdot \partial_s \wb^i \right)\, ds\\
 &+ \sum_{i=1}^{n_\cE} (\Pb^i_+ \cdot \vb^i(\ell^i) - \Pb^i_- \cdot \vb^i(0)) + \sum_{i=1}^{n_\cE}(\Qb^i_+ \cdot \wb^i(\ell^i) - \Qb^i_- \cdot \wb^i(0))\\
 &-\left(\Abb^+_\cI \Pb_+ -\Abb^-_\cI \Pb_- \right)\cdot \Vb - \left(\Abb^+_\cI \Qb_+ - \Abb^-_\cI \Qb_-\right)\cdot \Wb= - \sum_{i=1}^{n_\cE} \int_0^{\ell^i} \fb^i \cdot \vb^i\, ds,\\
 &\qquad \vb^i, \wb^i \in H^1(0,\ell^i;\ZR^3),\  i=1,\ldots, n_\cE, \quad \Vb, \Wb \in \ZR^{3n_\cV}.
\endaligned
\end{equation}
We then add (\ref{displ2}) and (\ref{rot2}) to (\ref{lj34}) and obtain
\begin{equation}\label{lj34A}
\aligned
 \sum_{i=1}^{n_\cE} & \int_0^{\ell^i} \Qbb^i (\Hbb^i)^{-1} (\Qbb^i)^T \qb^i \cdot \xib^i  -  (\partial_s \ub^i +\tb^i \times  \omegab^i) \cdot \thetab^i -\partial_s \omegab^i \cdot \xib^i\, ds \\
&+ \sum_{i=1}^{n_\cE} (\ub^i(\ell^i) \cdot \Thetab^i_+ - \ub^i(0) \cdot \Thetab^i_- ) + \sum_{i=1}^{n_\cE} (\omegab^i(\ell^i) \cdot \Xib^i_+ - \omegab^i(0) \cdot \Xib^i_- )\\
 & - ((\Abb^+_\cI)^T\Ub \cdot \Thetab_+ - (\Abb^-_\cI)^T\Ub \cdot \Thetab_-) - ((\Abb^+_\cI)^T\Omegab \cdot \Xib_+ - (\Abb^-_\cI)^T\Omegab \cdot \Xib_-)=0, \\
 &\qquad \xib^i, \thb^i \in  L^2(0,\ell^i;\ZR^3),\ i=1,\ldots, n_\cE,  \quad \Thetab_\pm, \Xib_\pm \in \ZR^{3n_\cE}.
\endaligned
\end{equation}
%\marginpar{novo}

In \cite{RadHAZU} and \cite{GT} the mixed formulation of the stent model was presented using the space $H^1(\cN;\ZR^3)$ for the displacement vector $\ub$ and the infinitesimal rotation vector $\omegab$. The space $L^2(\cN;\ZR^3)\times \ZR^3 \times \ZR^3$ for the Lagrange multipliers $\pb, \alb, \betab$, and  the continuity conditions for displacements and infinitesimal rotations were inherently built into the space $H^1(\cN;\ZR^3)$.  We now relax these conditions and consider them as additional equations in the problem and enlarge the space of unknowns by adding further Lagrange multipliers. The resulting function spaces are given by
\[
\aligned
V&=L^2(\cN;\ZR^3)\times L^2(\cN;\ZR^3)\times \ZR^{3n_\cE} \times \ZR^{3n_\cE} \times \ZR^{3n_\cE} \times \ZR^{3n_\cE} \times \ZR^{3} \times \ZR^{3}, \\
M&=L^2_{H^1}(\cN;\ZR^3) \times L^2_{H^1}(\cN;\ZR^3) \times \ZR^{3n_\cV} \times\ZR^{3n_\cV}.
\endaligned
\]
To simplify the notation for the elements of these spaces we introduce
\[
\Sib:=(\qb, \pb, \Pb_+,\Pb_-,\Qb_+,\Qb_-,\alb, \betab) \in V, \qquad \phib:=(\ub, \omegab, \Ub, \Omegab) \in M
\]
for the unknowns in the problem and
\[
\Gamb:=(\xib, \thetab, \Thetab_+,\Thetab_-,\Xib_+,\Xib_-,\gamb,\deltab) \in V, \qquad \psib:=(\vb, \wb, \Vb, \Wb) \in M
\]
for the associated test functions.

In this notation, the bilinear forms and the linear functionals that appear in the above calculations are given by
\[
\aligned
&a: V \times V \to \ZR, \qquad b : V \times M \to \ZR, \qquad f: M \to \ZR,\\
&a(\funS,\funG) :=  \sum_{i=1}^{n_\cE} \int_0^{\ell^i} \Qbb^i (\Hbb^i)^{-1} (\Qbb^i)^T \qb^i \cdot \xib^i\, ds,\\
&b(\funS,\funv)  := \sum_{i=1}^{n_\cE} \int_0^{\ell^i} \left(- \pb^i \cdot (\partial_s \vb^i +\tb^i \times \wb^i)  -  \qb^i  \cdot \partial_s \wb^i \right)\, ds\\
 &+ \sum_{i=1}^{n_\cE} (\Pb^i_+ \cdot \vb^i(\ell^i) - \Pb^i_- \cdot \vb^i(0)) + \sum_{i=1}^{n_\cE}(\Qb^i_+ \cdot \wb^i(\ell^i) - \Qb^i_- \cdot \wb^i(0))\\
 &-\left(\Abb^+_\cI \Pb_+ -\Abb^-_\cI \Pb_- \right)\cdot \Vb - \left(\Abb^+_\cI \Qb_+ - \Abb^-_\cI \Qb_-\right)\cdot \Wb + \alb \cdot \int_\cN \vb + \betab \cdot \int_\cN \wb,\\
 &f(\funv) := -\sum_{i=1}^{n_\cE} \int_0^{\ell^i} \fb^i \cdot \vb^i\, ds.
\endaligned
\]
Then the variational formulation (\ref{lj12A}), (\ref{lj34A}) and (\ref{*}) can be expressed as follows.

Determine $\funS \in V$ and $\funu \in M$ such that %\marginpar{Uvesti pokrate za nepoznanice?}
\begin{equation}\label{eqmain}
\aligned
a(\funS,\funG) + b(\funG,\funu) &= 0, & \funG \in V,\\
b(\funS,\funv) \phantom{a(\funS,\funG) + }&= f(\funv), & \funv \in M.
\endaligned
\end{equation}
In this way we have obtained that the solution of the stent problem as formulated in (\ref{lj1ri})--(\ref{fix}) satisfies (\ref{eqmain}) and conversely that any solution of (\ref{eqmain}) satisfies
(\ref{lj1ri})--(\ref{fix}).

In this section we have reformulated the mathematical formulation
of the stent model by including the continuity conditions at the nodes as extra equations and by adding further Lagrange multipliers. In the next sections, we will use this formulation to obtain a discrete $\inf$-$\sup$ inequality  and to present a simpler proof of the continuous $\inf$-$\sup$ inequality.

\section{Properties of the continuous model}
\setcounter{equation}{0}
\label{S4}
In this section we consider the properties of the continuous operator equation~\eqref{eqmain}. For the operator $B:V \to M'$ defined by
\[
\aligned
&{}_{M'}\langle B\funS, \funv\rangle_{M}=b(\funS,\funv), \qquad \funv \in M
\endaligned
\]
we have the adjoint operator $B^T: M \to V'$ (we use the matrix notation to illustrate the similarity to the discrete case discussed later), which satisfies
\[
\aligned
&{}_{V}\langle\funS, B^T\funv\rangle_{V'}=b(\funS,\funv), \qquad \funS \in V.
 \endaligned
\]
Then $\Ker{B^T}$, the kernel of $B^T$, is defined as a set of vector functions $\funv=(\vb,\wb,\Vb, \Wb)\in M$ such that
\[
b(\funS,\funv)=0, \qquad\funS \in V,
\]
so that $\funv=(\vb,\wb,\Vb, \Wb) \in \Ker{B^T}$ if and only if
\begin{eqnarray}\label{jdbe1a}
\partial_s \vb^i +\tb^i \times \wb^i &=&0, \quad \partial_s \wb^i =0,  \qquad i=1, \ldots, n_\cE,\\
\int_\cN \vb = \int_\cN \wb &=&0,\label{jdbe1b}\\
\vb^i(\ell^i) &=&\bbP^i_{\cE} (\Abb^+_\cI)^T \Vb, \quad \vb^i(0) =\bbP^i_{\cE} (\Abb^-_\cI)^T \Vb, \qquad i=1,\ldots,n_\cE,\label{jdbe1c}\\
\wb^i(\ell^i) &=&\bbP^i_{\cE} (\Abb^+_\cI)^T \Wb, \quad \wb^i(0) =\bbP^i_{\cE} (\Abb^-_\cI)^T \Wb, \qquad i=1,\ldots,n_\cE.\label{jdbe1d}
\end{eqnarray}
The  conditions \eqref{jdbe1c} and \eqref{jdbe1d} imply that $\vb$ and $\wb$ are continuous on the complete stent, i.e., $\vb,\wb\in H^1(\cN;\ZR^3)$. The conditions in (\ref{jdbe1a}) imply then that
%\begin{equation}\label{jdbe2}
%\partial_s \vb^i +\tb^i \times \wb^i =0, \quad \partial_s \wb^i =0,  \qquad i=1, \ldots, n_\cE.
%\end{equation}
%From the last equation
$\wb$ is constant on the complete stent and from  (\ref{jdbe1b})  we obtain $\wb=0$.
Analogously, from (\ref{jdbe1a}) we obtain that $\vb=0$, and hence $\Vb=\Wb=0$. Thus we have proved the following lemma.
\begin{lemma}\label{lKerBT}
$\Ker{B^T}=\{0\}$.
\end{lemma}
As next step we prove that $\Im{B^T}$ is closed. For this we derive a kind of Poincar\'{e} inequality on the graph $\cN$, using the notation $\vb(\ell) = [\vb^1(\ell^1), \ldots, \vb^{n_\cE}(\ell^{n_\cE})]^T$.
\begin{lemma}\label{ltheineq}
There exists a constant $C>0$ such that for all $\vb \in L^2_{H^1}(\cN;\ZR^3)$ and all $\Vb \in \ZR^{3n_\cV}$ the following inequality holds
\begin{equation}\label{poinin}
\aligned
\|\vb\|_{L^2(\cN;\ZR^3)} \leq C \bigg( &\|\vb'\|_{L^2(\cN;\ZR^3)}^2 +
\left( \int_\cN \vb\right)^2 + \left(\vb(0) - (\Abb^-_\cI)^T \Vb \right)^2\\
&+ \left(\vb(\ell) - (\Abb^+_\cI)^T \Vb \bigg)^2\right)^{1/2}.
\endaligned
\end{equation}
\end{lemma}
\proof
Suppose the contrary, i.e.,  for all constants $C>0$ there exist $\vb_C \in L^2_{H^1}(\cN;\ZR^3)$ and $\Vb_C \in \ZR^{3n_\cV}$ such that the opposite inequality holds. Then,  for $C=1/k$, $k\in \ZN$ there exist sequences $\vb_k \in L^2_{H^1}(\cN;\ZR^3)$ and $\Vb_k \in \ZR^{3n_\cV}$ such that
\[
\|\vb_k\|_{L^2(\cN;\ZR^3)} =1,
\]
\[
\|\vb_k'\|_{L^2(\cN;\ZR^3)}^2 +
\left( \int_\cN \vb_k\right)^2 + \left(\vb_k(0) - (\Abb^-_\cI)^T \Vb_k \right)^2+ \left(\vb_k(\ell) - (\Abb^+_\cI)^T \Vb_k \right)^2 \to 0,
\]
where, as before, $\vb_k'$ denotes the vector of partial derivatives of $\vb_k$.
Thus, taking an appropriate subsequence (still indexed by $k$), we have
\begin{eqnarray}\label{jdbe3a}
&&\vb_k \rightharpoonup \vb \mbox{ weakly in } L^2(\cN;\ZR^3),\\
&&\vb_k' \to 0 \mbox{ strongly in } L^2(\cN;\ZR^3),\label{jdbe3b}\\
&&\int_\cN \vb_k \to 0,\label{jdbe3c}\\
&&\vb_k(0) - (\Abb^-_\cI)^T \Vb_k \to 0,\label{jdbe3d}\\
&&\vb_k(\ell) - (\Abb^+_\cI)^T \Vb_k  \to 0.\label{jdbe3e}
\end{eqnarray}
It follows that on each strut we have
\[
\vb^i_k \rightharpoonup \vb^i \mbox{ weakly in } L^2(0,\ell^i;\ZR^3),\quad
\partial_s \vb^i_k \to 0 \mbox{ strongly in } L^2(0,\ell^i;\ZR^3), \quad i=1,\ldots,n_\cE,
\]
so that $\vb^i$ is constant on the $i$th strut and
\[
\vb^i_k \rightharpoonup \vb^i \mbox{ weakly in } H^1(0,\ell^i;\ZR^3).
\]
By the Trace Theorem; see e.g. \cite[Section 5.5, Theorem 1]{Evans} we have then
\[
\vb^i_k(0) \to \vb^i, \qquad \vb^i_k(\ell^i) \to \vb^i.
\]
Using (\ref{jdbe3d}) and (\ref{jdbe3e}), we have
\begin{equation}\label{**}
\bbP^i_\cE (\Abb^-_\cI)^T \Vb_k \to \vb^i,\quad
\bbP^i_\cE (\Abb^+_\cI)^T \Vb_k  \to \vb^i, \qquad i = 1,\ldots,n_\cE.
\end{equation}
Since in every block row of size $3$ the matrices $(\Abb^-_\cI)^T$ and $(\Abb^+_\cI)^T$ have exactly one identity matrix of size $3$ we have that $\Vb_k$ is convergent as well. We denote the limit by $\Vb$, and have that its values are given by $\vb^i$ suitably organized. Therefore, since $\Abb_\cI = \Abb^+_\cI-\Abb^-_\cI$, subtracting the sequences in (\ref{**}) we obtain that
\[
\Abb_\cI^T \Vb = 0.
\]
Since the rank of $\Abb_\cI^T$ is equal to $3n_\cV-3$,  the kernel of $\Abb_\cI^T$ is of dimension $3$ by the Rank--Nullity Theorem, \cite{Bapat}. We easily inspect that  $\Ker \Abb_{\cI}^T$ is spanned by the vectors
\[
(\eb_1, \eb_1, \ldots, \eb_1),\, (\eb_2, \eb_2, \ldots, \eb_2),\, (\eb_3, \eb_3, \ldots, \eb_3).
\]
Thus we obtain that all $\vb^i$ are equal. Since $0=\int_\cN \vb = \sum_{i=1}^{n_\cE} \ell^i \vb^i$, we obtain that
$\vb^i =0$ and hence $\Vb=0$, which also implies that $\vb=0$. Thus, since
\[
\vb^i_k(x) = \vb^i_k(0) + \int_0^{x} \partial_s \vb^i_k(s)\, ds,
\]
$(v^i_k)_k$ tends to 0 strongly in $L^2(0,\ell^i;\ZR^3)$ for all $i=1,\ldots,n_\cE$, which is in contradiction to the unit norm assumption of the sequence, i.e.,  $\|\vb_k\|_{L^2(\cN;\ZR^3)}=1$.
\proof

\begin{lemma}\label{ImBT}
$\Im{B^T}$ is closed.
\end{lemma}
\proof
Consider a convergent sequence in $\Im{B^T}$, i.e., a sequence of the form
\begin{equation}\label{jdbe5}
\aligned
&\partial_s \vb^i_k +\tb^i \times \wb^i_k\to \pb^i , \quad \partial_s \wb^i_k \to \qb^i,  \quad \mbox{strongly in } L^2(0,\ell^i;\ZR^3) \quad i=1, \ldots, n_\cE,\\
&\int_\cN \vb_k \to \alb,  \int_\cN \wb_k \to \betab,\\
&\vb^i_k(\ell^i) - \bbP^i_{\cE} (\Abb^+_\cI)^T \Vb_k \to \Pb^i_+, \qquad \vb^i_k(0) -\bbP^i_{\cE} (\Abb^-_\cI)^T \Vb_k \to \Pb^i_-, \\
&\wb^i_k(\ell^i) -\bbP^i_{\cE} (\Abb^+_\cI)^T \Wb_k \to \Qb^i_+, \qquad \wb^i_k(0) -\bbP^i_{\cE} (\Abb^-_\cI)^T \Wb_k \to \Qb^i_-.
\endaligned
\end{equation}
Then applying the inequality (\ref{poinin}) to the sequences $\wb_k=(\wb^1_{k},\ldots,\wb^{n_\cE}_k)$ and
$\Wb_k=(\Wb^1_k,\ldots,\Wb^{n_\cE}_k)$ implies that $\wb_k$ is bounded in $L^2(\cN,\ZR^3)$. Therefore, $\wb_k^i$ is bounded in $H^1(0,\ell^i;\ZR^3)$ and hence there exists a subsequence and a function $\wb^i\in H^1(0,\ell^i;\ZR^3)$ such that
\[
\wb^i_{k_l} \rightharpoonup \wb^i \mbox{ weakly in } H^1(0,\ell^i;\ZR^3), \qquad i=1,\ldots, n_\cE.
\]
We collect the limits in  $\wb=(\wb^1,\ldots,\wb^{n_\cE})$ and, using again the Trace Theorem, we obtain that
\[
\qb^i=\partial_s \wb^i, \quad  -\bbP^i_{\cE} (\Abb^+_\cI)^T \Wb_{k_l} \to \Qb^i_+ -\wb^i(\ell^i), \quad  -\bbP^i_{\cE} (\Abb^-_\cI)^T \Wb_{k_l} \to \Qb^i_- -\wb^i(0),
\]
for $i=1,\ldots,n_\cE$ and
\[
\betab=\int_\cN \wb.
\]
Since in each block row of dimension $3$ the matrices $\Abb^\pm_\cI$ have exactly one identity matrix, we obtain that $\Wb_{k_l}$ converges to $\Wb$ which satisfies
\[
-\bbP^i_{\cE} (\Abb^+_\cI)^T \Wb = \Qb^i_+ -\wb^i(\ell^i), \quad -\bbP^i_{\cE} (\Abb^-_\cI)^T \Wb = \Qb^i_- -\wb^i(0), \qquad i=1,\ldots, n_\cE.
\]
By inequality (\ref{poinin}), there is a unique function $\wb$ that satisfies the associated homogeneous system
\[
\partial_s \wb^i = 0, \ i=1,\ldots, n_\cE, \quad \int_\cN  \wb = 0, \quad \vb(0) - (\Abb^-_\cI)^T \Vb = \vb(\ell) - (\Abb^+_\cI)^T \Vb =0.
\]
Therefore, $\wb$ is unique and thus the whole sequences $(\wb_k)_k$ and $(\Wb_k)_k$ are convergent. Application of Lemma~\ref{ltheineq} to $\wb_k-\wb$ and $\Wb_k-\Wb$ implies that $\wb_k \to \wb$ strongly in $L^2_{H^1}(\cN;\ZR^3)$. Once we have this convergence we apply it to the term $\tb^i \times \wb^i_k$ and by the same reasoning identify all limits related to $\vb_k$ and $\Vb_k$.
We obtain $\vb_k \to \vb$ strongly in $L^2_{H^1}(\cN;\ZR^3)$ and $\Vb_k \to \Vb$ in $\ZR^{3n_\cV}$ and that
\[
\partial_s \vb^i +\tb^i \times \wb^i = \pb^i, \ \vb^i(\ell^i) - \bbP^i_{\cE} (\Abb^+_\cI)^T \Vb =  \Pb^i_+, \ \vb^i (0) -\bbP^i_{\cE} (\Abb^-_\cI)^T \Vb = \Pb^i_-, \ \alb=\int_\cN \ub.
\]
Thus $\funS=(\qb,\pb,\Pb_+,\Pb_-,\Qb_+,\Qb_-,\alb,\betab)$ belongs to $\Im{B^T}$,
and hence  $\Im{B^T}$ is closed.
\eproof

As a direct consequence of Lemmas~\ref{lKerBT},~\ref{ImBT}, and \cite[Proposition 1.2, page 39]{BrezziFortin}, we obtain the following corollary.
\begin{corollary}[Continuous $\inf$-$\sup$ inequality]\label{cinfsup}
Consider the variational formulation of the stent model~\eqref{eqmain}. Then there exists a constant $k_c>0$ such that
\[
\inf_{\funv\in M} \sup_{\funS\in V} \frac{b(\funS,\funv)}{\|\funS\|_V \|\funv\|_M} \geq k_c.
\]
\end{corollary}
As our next result we will prove the $\Ker{B}$ ellipticity of the form $a$, i.e., that there exist $c_a>0$ such that
\[
a(\funS,\funS) \geq c_a \|\funS\|_{V}^2, \qquad \funS \in \Ker{B}.
\]
To obtain this result, we need to restrict the class of networks we consider.

\begin{lemma}\label{lelliptic}
Let the stent geometry be  such that
\[
\sum_{i\in J^+_j \atop  i \text{th edge is straight}} \alpha_i \tb^i - \sum_{i\in J^-_j \atop i \text{th edge is straight}} \alpha_i \tb^i = 0, \quad j=1,\ldots, n_\cV
\]
implies that $\alpha_i=0$ for all straight edges $i$.
Then the bilinear form  $a$ from the variational formulation \eqref{eqmain} is $\Ker{B}$ elliptic.
\end{lemma}
\proof
The space $\Ker{B}$ is given by the set of $\funS=(\qb, \pb, \Pb_+,\Pb_-,\Qb_+,\Qb_-,\alb, \betab)\in V$ such that
\[
b(\funS,\funv) = 0, \qquad \funv\in M.
\]
Thus, from the definition of the form $b$ one has that for all $\funv=(\vb,\wb,\Vb,\Wb)\in M$
\begin{equation}\label{eKerB}
\aligned
&\sum_{i=1}^{n_\cE} \int_0^{\ell^i} \left(- \pb^i \cdot (\partial_s \vb^i +\tb^i \times \wb^i)  -  \qb^i  \cdot \partial_s \wb^i \right)\, ds\\
 &+ \sum_{i=1}^{n_\cE} (\Pb^i_+ \cdot \vb^i(\ell^i) - \Pb^i_- \cdot \vb^i(0)) + \sum_{i=1}^{n_\cE}(\Qb^i_+ \cdot \wb^i(\ell^i) - \Qb^i_- \cdot \wb^i(0))\\
 &-\left(\Abb^+_\cI \Pb_+ -\Abb^-_\cI \Pb_- \right)\cdot \Vb - \left(\Abb^+_\cI \Qb_+ - \Abb^-_\cI \Qb_-\right)\cdot \Wb + \alb \cdot \int_\cN \vb + \betab \cdot \int_\cN \wb=0.
 \endaligned
\end{equation}
Since $\Vb$ and $\Wb$ in $\ZR^{n_\cV}$ are arbitrary, we obtain
\begin{equation}\label{eeq2}
\Abb^+_\cI \Pb_+ -\Abb^-_\cI \Pb_- = \Abb^+_\cI \Qb_+ - \Abb^-_\cI \Qb_-=0,
\end{equation}
which is equivalent to (\ref{sila}).
These conditions mean that at each vertex the sum of the contact forces as well as the sum of the contact moments are zero. % (dynamic junction conditions).
Then, for $\gamb \in \ZR^3$, we insert $\vb^i=\gamb$, $\wb^i=0$, $\Wb=\Vb=0$ as test functions in (\ref{eKerB}) and obtain
\[
\sum_{i=1}^{n_\cE} (\Pb^i_+  - \Pb^i_-) \cdot \gamb + \alb \cdot \left(\sum_{i=1}^{n_\cE}\ell^i \right)\gamb =0.
\]
Since from (\ref{eeq2})
\[
\sum_{i=1}^{n_\cE} (\Pb^i_+  - \Pb^i_-) = \sum_{i=1}^{n_\cE} \Pbb^i_\cE (\Abb^+_\cI \Pb_+ -\Abb^-_\cI \Pb_-) =0,
\]
we obtain that $\alb =0$. Now, for fixed $i$ we insert {$\vb^i \in C^1([0, \ell^i];\ZR^3)$ with compact support in $\langle0,\ell^i\rangle$} in (\ref{eKerB})  and obtain $\int_0^{\ell^i} \pb^i \cdot \partial_s \vb^i=0$. This implies that $\pb^i$ is a constant on each strut. Inserting a single {$\vb^i \in C^1([0,\ell^i];\ZR^3)$} in (\ref{eKerB}), we obtain that
\begin{equation}\label{***}
-\pb^i \cdot  \int_0^{\ell^i} \partial_s \vb^i \, ds + \Pb^i_+ \cdot \vb^i(\ell^i) - \Pb^i_- \cdot \vb^i(0) =0
\end{equation}
and thus $\pb^i = \Pb^i_+ = \Pb^i_-$.

Similarly, for $\gamb\in \ZR^3$ we insert $\wb^i=\gamb$ for all $i$ in (\ref{eKerB}) and obtain
\[
%\begin{equation}\label{eKerB}
- \sum_{i=1}^{n_\cE} \int_0^{\ell^i}  \pb^i \times \tb^i  \cdot\gamb \, ds + \sum_{i=1}^{n_\cE}(\Qb^i_+ - \Qb^i_-) \cdot \gamb + \betab \cdot \left(\sum_{i=1}^{n_\cE}\ell^i \right) \gamb=0.
\]
%\end{equation}
As in the case of contact forces we have $\sum_{i=1}^{n_\cE}(\Qb^i_+ - \Qb^i_-) =0$. For the first term we argue as follows
\[
\sum_{i=1}^{n_\cE} \int_0^{\ell^i}  \pb^i \times \tb^i\, ds = \sum_{i=1}^{n_\cE} \pb^i  \times  (\Phib^i(\ell^i) -\Phib^i(0)) = \sum_{j=1}^{n_\cV} (\sum_{i\in J_j^+}\Pb^i_+ - \sum_{i\in J^-_j}\Pb^i_-) \times \cV_j=0,
\]
again by (\ref{eeq2}), where $\cV_j$ denotes the $j$th vertex. Thus we conclude that $\betab=0$. What is left from (\ref{eKerB}) are the equations for all $i\in \{1, \ldots, n_\cE\}$ and all $\wb^i\in H^1(0,\ell^i)$ given by
\begin{equation}\label{eKerB1}
- \int_0^{\ell^i} \pb^i \times \tb^i \cdot  \wb^i\, ds   -  \int_0^{\ell^i} \qb^i  \cdot \partial_s \wb^i \, ds + \Qb^i_+ \cdot \wb^i(\ell^i) - \Qb^i_- \cdot \wb^i(0)=0.
\end{equation}
Defining $\qtb^i = \qb^i- \pb^i \times (\Phib^i(s)-\Phib^i(0))$ and inserting this form in (\ref{eKerB1}), we obtain
\[%\begin{equation}\label{eKerB2}
\aligned
&- \int_0^{\ell^i}\pb^i \times \tb^i \cdot  \wb^i\, ds   -  \int_0^{\ell^i} \qtb^i  \cdot \partial_s \wb^i\, ds  -\int_0^{\ell^i} \pb^i \times (\Phib^i(s)-\Phib^i(0)) \cdot \partial_s \wb^i\, ds\\
&\qquad + \Qb^i_+ \cdot \wb^i(\ell^i) - \Qb^i_- \cdot \wb^i(0)=0.
\endaligned
\]
%\end{equation}
After partial integration in the third term we obtain
\begin{equation}\label{eKerB2}
-  \int_0^{\ell^i} \qtb^i  \cdot \partial_s \wb^i \, ds + (\Qb^i_+ - \pb^i \times (\Phib^i(\ell^i)-\Phib^i(0)))\cdot \wb^i(\ell^i) - \Qb^i_- \cdot \wb^i(0)=0.
\end{equation}
As in the case of contact forces, see \eqref{***}, this implies that the $\qtb^i$ are constants and that
\[
\qtb^i =\Qb^i_+ -\pb^i \times (\Phib^i(\ell^i)-\Phib^i(0)) = \Qb^i_-, \qquad i=1,\ldots, n_\cE.
\]
Thus, we have obtained the following characterization of $\Ker{B}$,
\[
\aligned
&\pb^i = \Pb^i_+ = \Pb^i_-, \qquad  \qb^i = \pb^i \times (\Phib^i(s)-\Phib^i(0)) + \Qb^i_-, \\
&\Qb^i_+ = \pb^i \times (\Phib^i(\ell^i)-\Phib^i(0))  + \Qb^i_-, \qquad \alb=\betab=0,
\endaligned
\]
with (\ref{eeq2}) being  satisfied.

Since $\Ker{B}$ is of finite dimension and the form $a$ is obviously positive semidefinite, to check that $a$ is $\Ker{B}$ elliptic, we only have to check that $\Ker{a}\cap \Ker{B}$ is trivial.
Assume that $\funS \in \Ker{a} \cap \Ker{B}$. For elements of $\Ker{a}$, $\qb^i=0, i=1,\ldots,n_\cE$  and for elements in $\Ker{B}$ we further have
\[
\pb^i \times (\Phib^i({s})-\Phib^i(0)) =  \Qb^i_- = 0 \quad s \in [0,\ell^i], \qquad \Qb^i_+=\alb=\betab=0, \qquad i=1,\ldots, n_\cE.
\]
From the first equation we have that $\pb^i = \alpha_i \tb^i$, for some $\alpha_i \in \ZR$, $i=1,\ldots,n_\cE$ if the strut is straight, otherwise $\pb^i=0$. By our assumption on the geometry of the stent this implies that  $\Ker{a}\cap \Ker{B}=\{0\}$. This concludes the proof.
\eproof

The restriction on the geometry in Lemma~\ref{cinfsuph} does not exclude typical examples of stents, since most of the struts in stents are curved. But even if they were straight, then we can start from the vertices where only two struts meet and conclude that the associated  scalars $\alpha_i$ for these struts should be zero and then continue until we conclude that all coefficients are zero.

The properties of the continuous model that we have shown imply {the existence and} the uniqueness of the solution.
\begin{theorem}\label{main}
There is a unique solution of (\ref{eqmain}).
\end{theorem}
\proof
Since $a$ is $\Ker{B}$ elliptic by Lemma~\ref{lelliptic} and {since} the $\inf-\sup$ inequality holds on $V\times M$ by Lemma~\ref{lKerBT} and Lemma~\ref{ImBT},  the application of \cite[Corollary 4.1, page 61]{GR} or \cite[Theorem 1.1, page 42]{BrezziFortin} implies the assertion.
\eproof

Note that, even though we deal with the pure traction problem, because of the introduction of two additional conditions on total displacement and total infinitesimal rotation in \eqref{fix}, we have a unique solution of (\ref{eqmain}) for all forces. Furthermore, the necessary conditions usual for the pure traction problem (zero total force and zero total couple) are not necessary any more. The Lagrange multipliers $\alb$ and $\betab$ deal with that. See \cite{RadHAZU} for explicit formulas for these multipliers.

\section{Properties of the discrete model}
\setcounter{equation}{0}
\label{S5}
In this section we discuss a discrete approximation of the problem (\ref{eqmain}). The derivations in this section are done for straight struts, i.e., under the assumption that $\ell^i \tb^i =\Phib^i(\ell^i)-\Phib^i(0)$ for all struts.  %$i$'s.
If a strut is curved then  it is approximated by a piecewise straight approximation. In this case we need estimates for two solutions of the stent model for two geometries.

Let us denote by $P_m(\cN)$ the set of functions on the graph $\cN$ which are polynomials of degree $m \geq 0$ on each edge of the graph. Note that we do not assume that the functions in $P_m(\cN)$ are continuous at vertices. We use the similar notation $P_m([0,\ell])$ for the space of polynomials of degree $m$ on the segment $[0,\ell]$. For $k,n\in \ZN_{0}$ we define the finite dimensional spaces of discrete approximations
\[
\aligned
&V_k := P_k(\cN)^3 \times P_k(\cN)^3 \times \ZR^{3n_\cE} \times \ZR^{3n_\cE} \times \ZR^{3n_\cE}\times \ZR^{3n_\cE}\times \ZR^{3} \times \ZR^{3}, \\
&M_n := P_n(\cN)^3 \times P_n(\cN)^3 \times \ZR^{3n_\cV} \times \ZR^{3n_\cV}.
\endaligned
\]
%
%Next we choose subspaces $V_k\subset V$ and $M_n\subset M$ consisting of polynomials of degree $k$ %on each strut for %the components $\pb$ and $\qb$ in $V_k$ and polynomials of degree $n$ on each %strut for $\vb$ and $\wb$ in $M_n$.
%These spaces are finite-dimensional.
We assume that $n\geq 1$ and $k\geq 0$ and consider the following discrete approximation of (\ref{eqmain}).

Determine $\funS\in V_k$ and $\funu \in M_n$ such that
\begin{equation}\label{eqmainh}
\aligned
a(\funS,\funG) + b(\funG,\funu) &= 0, \qquad \qquad\funG \in V_k,\\
b(\funS,\funv) \phantom {+a(\funS,\funG)}&= f(\funv), \qquad \funv \in M_n.
\endaligned
\end{equation}
%
%In this section we denote by $P_m$ the set of all polynomials of degree less then $m$.
%
The form $b$ on $V_k\times M_n$ defines the operator $B_h:V_k\to M_n'$, where $M_n'$ denotes the dual of $M_n$. In general, however, $\Ker{B_h}$ is not a subset of $\Ker{B}$. However, we will show that  if $n-1\geq k$ then  it is  a subset and thus applying Lemma~\ref{lelliptic} gives the following result.
\begin{lemma}\label{lelliptich}
Consider the discrete problem \eqref{eqmainh} and let $n-1\geq k$. Then the bilinear form $a$ is $\Ker{B_h}$ elliptic.
\end{lemma}
\proof
As in the continuous case, the elements $\funS=(\qb, \pb, \Pb_+,\Pb_-,\Qb_+,\Qb_-,\alb, \betab)$ of $\Ker{B_h}$ satisfy (\ref{eeq2})
%\begin{equation}\label{eeq2h}
%\Abb^+_\cI \Pb_+ -\Abb^-_\cI \Pb_- = \Abb^+_\cI \Qb_+ - \Abb^-_\cI \Qb_-=0.
%\end{equation}
%As in the proof of Lemma~\ref{lelliptic} we first obtain (\ref{eeq2})
and by the same arguments it follows that $\alb=0$. For fixed $i$ and a test function $\vb^i \in P_n([0,\ell^i])$ we obtain the equation
%\begin{equation}\label{eKerB}
\[
-\int_0^{\ell^i}  \pb^i \cdot \partial_s \vb^i \, ds  + \Pb^i_+ \cdot \vb^i(\ell^i) - \Pb^i_- \cdot \vb^i(0) =0.
\]
%\end{equation}
For constant $\vb^i = \gamb \in \ZR^3$, we obtain $\Pb^i_+=\Pb^i_-$ and,
inserting $\pb^i=\ptb^i+\Pb^i_-$ implies
\[%\begin{equation}\label{eKerB}
-\int_0^{\ell^i}  \ptb^i \cdot \partial_s \vb^i\, ds =0,\quad \vb^i \in P_n([0,\ell^i]).
\]
%\end{equation}
Since $n\geq k+1$, the function $\ptb^i$ is zero and hence $\pb^i=\Pb^i_+=\Pb^i_-$.

For $\vb^i=0$, $i=1,\ldots,n_\cE$ and $\Vb=\Wb=0$ in (\ref{eKerB}) we obtain
\begin{equation}\label{eKerBh}
\sum_{i=1}^{n_\cE} \int_0^{\ell^i} - \pb^i \cdot \tb^i \times \wb^i  -  \qb^i  \cdot \partial_s \wb^i \, ds + \sum_{i=1}^{n_\cE}(\Qb^i_+ \cdot \wb^i(\ell^i) - \Qb^i_- \cdot \wb^i(0))+ \betab \cdot \int_\cN \wb=0.
\end{equation}
Inserting $\wb^i=\gamb \in \ZR^3$, $i=1,\ldots,n_\cE$, we obtain
\[
\sum_{i=1}^{n_\cE} -  \gamb \cdot  \int_0^{\ell^i} \pb^i \times \tb^i \, ds  + \sum_{i=1}^{n_\cE}(\Qb^i_+  - \Qb^i_- )\cdot \gamb+ \betab \cdot \left( \sum_{i=1}^{n_\cE} \ell^i\right)\gamb=0.
\]
As in the proof of Lemma~\ref{lelliptic}, we obtain that $\betab=0$, and thus we are left with the equation
\begin{equation}\label{eKerBh1}
- \pb^i \times \tb^i   \cdot \int_0^{\ell^i} \wb^i \, ds -  \int_0^{\ell^i} \qb^i  \cdot \partial_s \wb^i \, ds+ \Qb^i_+ \cdot \wb^i(\ell^i) - \Qb^i_- \cdot \wb^i(0)=0, \qquad \wb^i \in P_n([0,\ell^i]).
\end{equation}
For $k\geq 1$, we insert  $\qb^i=\qtb^i+s \pb^i \times \tb^i$ into this equation and obtain
\[
%\begin{equation}\label{eKerBh1}
- \pb^i \times \tb^i   \cdot \int_0^{\ell^i} \wb^i\, ds -  \int_0^{\ell^i} \qtb^i \cdot \partial_s \wb^i\, ds  - \pb^i \times \tb^i   \cdot \int_0^{\ell^i} s  \partial_s \wb^i \, ds + \Qb^i_+ \cdot \wb^i(\ell^i) - \Qb^i_- \cdot \wb^i(0)=0.
\]
%\end{equation}
After partial integration in the third term we obtain
\[
%\begin{equation}\label{eKerBh1}
-  \int_0^{\ell^i} \qtb^i \cdot \partial_s \wb^i\,ds  + (\Qb^i_+ - \ell^i  \pb^i \times \tb^i )\cdot \wb^i(\ell^i) - \Qb^i_- \cdot \wb^i(0)=0.
\]
%\end{equation}
Since constant functions  are contained in $V_k$, for $\wb^i = \gamb$ we obtain $\Qb^i_+ = \Qb^i_- + \ell^i  \pb^i \times \tb^i$. Then setting $\qtb^i = \vec{\tilde{\tilde{q}}} + \Qb^i_-$, we obtain
\[
\int_0^{\ell^i} \vec{\tilde{\tilde{q}}}^i \cdot \partial_s \wb^i\, ds =0, \qquad \wb^i \in P_n([0,\ell^i]).
\]
As before, since $n -1\geq k$, this implies that $\vec{\tilde{\tilde{q}}}=0$, and hence $\qb^i=\Qb^i_- +s \pb^i \times \tb^i$.

For $k=0$, $\qb^i$ is constant, so from (\ref{eKerBh1}) we obtain
\[
%\begin{equation}\label{eKerBh1}
- \pb^i \times \tb^i   \cdot \int_0^{\ell^i} \wb^i \, ds - \qb^i  \cdot  (\wb^i(\ell^i)-\wb^i(0)) + \Qb^i_+ \cdot \wb^i(\ell^i) - \Qb^i_- \cdot \wb^i(0)=0.
\]
%\end{equation}
This implies that
\[
\pb^i \times \tb^i=0, \qquad \qb^i = \Qb^i_+=\Qb^i_-
\]
and we have obtained the characterization of $\Ker{B_h}$ given by $(\qb,\pb,\Pb_+,\Pb_-,\Qb_+,\Qb_-,\alb,\betab)$ that satisfy (\ref{eeq2}) and
\[
\pb^i=\Pb^i_+=\Pb^i_-, \qquad \alb=\betab=0, \qquad \qb^i = \Qb^i_- +s \pb^i \times \tb^i, \qquad \Qb^i_+=\Qb^i_- + \ell^i  \pb^i \times \tb^i.
\]
Additionally, if $k=0$, then $\pb^i \times \tb^i=0$.
Thus, $\Ker{B_h} = \Ker{B} \cap V_k \subseteq \Ker{B}$ and hence $a$ is elliptic on $\Ker{B_h}$ by Lemma~\ref{lelliptic}.
\eproof

\begin{lemma}\label{lKerBTh}
Consider the discrete problem \eqref{eqmainh} and let $k\geq n-1$. Then $\Ker{B_h^T}=\{0\}$.
\end{lemma}
\proof
$\Ker{B_h^T}$ is defined as a set of $\funv=(\vb,\wb,\Vb, \Wb)\in M_n$ such that
\[
b(\funS,\funv)=0, \qquad \funS =(\qb, \pb, \Pb_+,\Pb_-,\Qb_+,\Qb_-,\alb, \betab) \in V_k.
\]
Thus, $\funv=(\vb,\wb,\Vb, \Wb) \in \Ker{B_h^T}$ if and only if
\[
\aligned
&\hspace{-2ex}\sum_{i=1}^{n_\cE} \int_0^{\ell^i} - \pb^i \cdot (\partial_s \vb^i +\tb^i \times \wb^i)  -  \qb^i  \cdot \partial_s \wb^i \, ds\\
 &+ \sum_{i=1}^{n_\cE} (\Pb^i_+ \cdot \vb^i(\ell^i) - \Pb^i_- \cdot \vb^i(0)) + \sum_{i=1}^{n_\cE}(\Qb^i_+ \cdot \wb^i(\ell^i) - \Qb^i_- \cdot \wb^i(0))\\
 &-\left(\Abb^+_\cI \Pb_+ -\Abb^-_\cI \Pb_- \right)\cdot \Vb - \left(\Abb^+_\cI \Qb_+ - \Abb^-_\cI \Qb_-\right)\cdot \Wb + \alb \cdot \int_\cN \vb + \betab \cdot \int_\cN \wb =0
\endaligned
\]
for all $\funS \in V_k$. This is equivalent to
\begin{eqnarray}\label{kerbhTa}
&&\sum_{i=1}^{n_\cE} \int_0^{\ell^i} \left(- \pb^i \cdot (\partial_s \vb^i +\tb^i \times \wb^i)  -  \qb^i  \cdot \partial_s \wb^i \right)\, ds =0, \pb^i, \qb^i \in P_k([0,\ell^i]),\\
&&\qquad   i=1, \ldots, n_\cE,\\
&&\int_\cN \vb = \int_\cN \wb =0,\label{kerbhTb}\\
&&\vb^i(\ell^i) =\bbP^i_{\cE} (\Abb^+_\cI)^T \Vb, \qquad \vb^i(0) =\bbP^i_{\cE} (\Abb^-_\cI)^T \Vb, \qquad i=1,\ldots,n_\cE,\label{kerbhTc}\\
&&\wb^i(\ell^i) =\bbP^i_{\cE} (\Abb^+_\cI)^T \Wb, \qquad \wb^i(0) =\bbP^i_{\cE} (\Abb^-_\cI)^T \Wb, \qquad i=1,\ldots,n_\cE.\label{kerbhTd}
\end{eqnarray}
Since $k\geq n-1$, from \eqref{kerbhTa} for a test function $\qb^i$ we obtain that $\wb^i$ is constant for each strut.
The continuity of the infinitesimal rotations at the vertices follows from \eqref{kerbhTd}. This
implies that $\wb^i=\const$ and then (\ref{kerbhTb}) implies that $\wb^i=0$, i.e., $\wb^i=0$ for all $i=1,\ldots,n_\cE$ and $\Wb=0$. Analogous arguments using \eqref{kerbhTc} imply that $\vb^i=0$ and $\Vb=0$ as well.
\eproof

Let $n=k+1$, so that both Lemma~\ref{lelliptich} and Lemma~\ref{lKerBTh} apply.
Since we are in the finite dimensional case, clearly $\Im{B_h}$ is closed. Proposition 1.2, page 39 in \cite{BrezziFortin} then implies that $\Im{B_h} = (\Ker{B_h^T})^0$ and then by Lemma~\ref{lKerBTh} it follows that $\Im{B_h} = M_n'$.
Furthermore, by Lemma~\ref{lelliptich}, the bilinear form  $a$ is $\Ker{B_h}$ elliptic. Then, by the classical theory for finite dimensional approximations of mixed formulations, e.g. Proposition 2.1 in \cite{BrezziFortin}, we obtain the following existence and uniqueness result for the discretized problem.
\begin{theorem}\label{texistenceh}
Let $n=k+1$. Then problem (\ref{eqmainh}) has a unique solution.
\end{theorem}
%
%\section{Error estimates}
%\setcounter{equation}{0}
%Next we prove discrete $\inf$-$\sup$ inequality.
%
Applying the classical results then also we obtain the discrete $\inf$-$\sup$ inequality.
\begin{corollary}[Discrete $\inf$-$\sup$ inequality]\label{cinfsuph}
If  $n=k+1$, then there exists a constant $k_d>0$ such that
\[
\inf_{\funv\in M_n} \sup_{\funS\in V_k} \frac{b(\funS,\funv)}{\|\funS\|_{V_k} \|\funv\|_{M_n}} \geq k_d.
\]
\end{corollary}
\proof
By Corollary~\ref{cinfsup}, the continuous $\inf$-$\sup$ inequality holds.
By Lemma~\ref{lKerBT} and Lemma~\ref{lKerBTh}  we have $\Ker{B_h^T} = \Ker{B^T}=\{0\}$. Thus Proposition 2.2, page 53 in \cite{BrezziFortin} implies that the assumptions of Proposition 2.8, page 58 in \cite{BrezziFortin} are fulfilled and we obtain the discrete $\inf$-$\sup$ inequality.
\eproof

\begin{remark}{\rm
Note that the constant $k_d$ from Corollary~\ref{cinfsuph} depends on the subspaces $V_k$ and $M_n$.
}
\end{remark}

Using Theorem 2.1, page 60 in \cite{BrezziFortin}, the discrete $\inf$-$\sup$ inequality in Corollary~\ref{cinfsuph} and Lemma~\ref{lelliptich}, i.e., the coercivity of the form $a$ on $\Ker{B_h}$, we obtain error estimates also in the discrete problem. Introducing analogous notation as in the continuous case,
\[
\Sib^h:=(\qb^h, \pb^h, \Pb^h_+,\Pb^h_-,\Qb^h_+,\Qb^h_-,\alb^h, \betab^h) \in V_k, \qquad \phib^h:=(\ub^h, \omegab^h, \Ub^h, \Omegab^h) \in M_n
\]
for the unknowns in the problem and
\[
\Gamb^h:=(\xib^h, \thetab^h, \Thetab^h_+,\Thetab^h_-,\Xib^h_+,\Xib^h_-,\gamb^h,\deltab^h) \in V_k, \qquad \psib^h:=(\vb^h, \wb^h, \Vb^h, \Wb^h) \in M_n
\]
for the test functions we have the following theorem.
\begin{theorem}\label{terror}
Let $n=k+1$ and let $(\funS, \funu)\in V\times M$ be the solution of (\ref{eqmain}) and let $(\funS^h, \funu^h)\in V_k\times M_n$ be the solution of (\ref{eqmainh}). Then
\[
\aligned
&\|\funS- \funS^h\|_{V} + \|\funu-\funu^h\|_{M} \leq c \Big( \inf_{\funG^h\in V_k}\|\funS - \funG^h\|_V+\inf_{\funv^h \in M_n}\|\funu - \funv^h\|_{M}\Big).
\endaligned
\]
\end{theorem}

\begin{remark}\label{rfinoca}\rm
The construction of finite elements, as presented, is directly related to the struts which are described
by their prescribed length. To increase the accuracy we can increase the polynomial degree, assuming that the geometry has been described without error. On the other hand, we can
change the topology of the stent by adding new points on existing struts (and thus not changing the geometry of the stent) in the original definition of the graph $\cN=(\cV,\cE)$.
In this way we obtain a refined model with transmission conditions of continuity of displacements, rotations, contact moments and forces at new points. Since at each new vertex only two struts meet, these transmission conditions are the same coupling conditions (kinematical and
dynamical) as for all vertices of the stent. Thus, the resulting
weak formulations as in \cite{IMASJ} or \cite{RadHAZU} are the same for both networks, the original one and the one with added vertices.
%For this network with added points we apply the discretization with prescribed mesh size $h$.
\end{remark}

The error estimate in Theorem~\ref{terror} can be employed in the finite element method using interpolation estimates in $L^2$ and $H^1$. Note that in contrast to the numerical method in \cite{GT}, due to the availability of the discrete $\inf$-$\sup$ inequality in the new formulation, we obtain the error estimate for all variables, including the contact forces.
It is a classical result, see e.g. \cite{CiarletFEM}, that for a function $\vphi\in H^r(0,\ell)$ and its polynomial Lagrange interpolant $\Pi_m\vphi$ of degree $m$ one has the estimate
\begin{equation}\label{pint}
\aligned
\|\vphi - \Pi_m\vphi\|_{L^2(0,\ell)} \leq C \ell^{\min{\{r,m+1\}}} \|\vphi^{(\min{\{r,m+1\}})}\|_{L^2(0,\ell)}, \qquad \vphi \in H^{r}(0,\ell),\\
\|\vphi - \Pi_m\vphi\|_{H^1(0,\ell)} \leq C \ell^{\min{\{r-1,m\}}} \|\vphi^{(\min{\{r,m+1\}})}\|_{L^2(0,\ell)}, \qquad \vphi \in H^{r}(0,\ell).\\
\endaligned
\end{equation}
With $ h:=\max\{\ell^i, i=1,\ldots,n_\cE\}$, then combining \eqref{pint} with Theorem~\ref{terror},
we obtain the following error estimate for the finite element method.
\begin{theorem}\label{terrorglavni}
Let $n=k+1$, $r\geq 0$, $r\in \ZN$, and let $\fb \in L^2_{H^r}(\cN;\ZR^3)$. Let $(\funS, \funu)\in V\times M$ be the solution of (\ref{eqmain}) and $(\funS^h, \funu^h)\in V_k\times M_n$ the solution of (\ref{eqmainh}). Then
\[
\aligned
\|\funS- \funS^h\|_{V} + \|\funu-\funu^h\|_{M}
&\leq c  h^{\min{\{r+1,k+1\}}} \|\fb\|_{L^2_{H^{r}}}.
\endaligned
\]
\end{theorem}
\proof
The error of the finite element approximation is estimated by the error of the interpolation operator. Thus we get
\[
\aligned
\|\funS- \funS^h\|_{V} + \|\funu-\funu^h\|_{M} &\leq c \Big( \|\funS - \Pi_k\funS\|_V+\|\funu - \Pi_n\funu\|_{M}\Big).
\endaligned
\]
Since in this section the struts are assumed to be straight for $\fb \in L^2_{H^r}(\cN;\ZR^3)$, from the differential equations we obtain that
\[
\pb \in L^2_{H^{r+1}}(\cN;\ZR^3), \quad \qb \in L^2_{H^{r+2}}(\cN;\ZR^3), \quad \omegab \in L^2_{H^{r+3}}(\cN;\ZR^3), \quad \ub \in L^2_{H^{r+4}}(\cN;\ZR^3).
\]
Thus, we obtain the estimate
\[
\aligned
&\hspace{-1ex}\|\funS- \funS^h\|_{V} + \|\funu-\funu^h\|_{M} \\
&\leq c \Big( h^{\min{\{r+1,k+1\}}}\|{(\qb,\pb)}^{(\min{\{r+1,k+1\}})}\|_{L^2(\cN;\ZR^6)}\\
 &\qquad \quad + h^{\min{\{r+2,n\}}}\|{(\ub,\omegab)}^{(\min{\{r+3,n+1\}})}\|_{L^2(\cN;\ZR^6)}\Big)\\
&\leq c \Big( h^{\min{\{r+1,k+1\}}}\|\fb\|_{L^2_{H^{\min{\{r,k\}}}}} + h^{\min{\{r+2,n\}}}\|\fb
\|_{L^2_{H^{\min{\{r,n-2\}}}}}\Big)\\
&\leq c  h^{\min{\{r+1,k+1\}}} \|\fb\|_{L^2_{H^{r}}}. 
\endaligned
\]
\eproof

Note that for $k=1$ and $n=2$ and $r\geq 1$ we obtain quadratic convergence in the $H^1$ norm for the displacements and infinitesimal rotations. This is in accordance with the convergence rate obtained in \cite{GT} for the classical formulation.

Having obtained error estimates for the continuous and discrete problem, in the next subsection, we
consider the properties of the resulting linear system.

\subsection{Block structure of the discretization matrix}
In the sequel we assume that $n=k+1$.
Using the same structure as in the continuous problem, the discrete problem is given by a linear system $\Kbb \xb = \Fb$, where
\[
\Kbb = \left[
\begin{array}{cc}
\Abb & \Bbbb^T\\
\Bbbb & 0
\end{array}
\right],\qquad \Fb=\left[\begin{array}{c}
0\\
\Fb_2
\end{array}\right]
\]
with a square matrix $\Abb$ of size $3(2k+6)n_\cE + 6$ and a rectangular matrix $\Bbbb$  of size $(3(2k+4)n_\cE+6 n_\cV)\times (3(k+6)n_\cE + 6)$. Having in mind the evolution problem that we will study in the next section, we partition these matrices further as
\begin{equation}\label{blstruc}
\aligned
\Abb = \left[
\begin{array}{cc}
\Abb_{11}& 0 \\
0 & 0
\end{array}
\right], \qquad \Bbbb= \left[\begin{array}{cc}
0 & \Bbbb_{32}\\
\Bbbb_{41} & \Bbbb_{42}
\end{array}\right],
\endaligned
\end{equation}
where
$\Abb_{11}$ is a square matrix of size $3(k+1)n_\cE$, $\Bbbb_{32}$ is a matrix of size $3(k+2)n_\cE \times (3(k+5)n_\cE +6)$, $\Bbbb_{41}$ is of size $(3(k+2)n_\cE+6n_\cV ) \times 3(k+1)n_\cE$ and $\Bbbb_{42}$ is of size $(3(k+2)n_\cE+6n_\cV) \times (3(k+5)n_\cE +6)$ associated with the following variables,
\begin{center}
\begin{tabular}{c||c||c||c||c||c|}
dim $\backslash$ unknown&$\qb$ & $(\pb,\Pb_+,\Pb_-,\Qb_+,\Qb_-,\alb,\betab)$ & $\ub$ & $(\omegab,\Ub,\Omegab)$ & $\Fb$\\\hline\hline
$3(k+1)n_\cE$ &$\Abb_{11}$ & 0 & 0 & $\Bbbb_{41}^T$& 0\\\hline \hline
$3(k+5)n_\cE+6$ & 0 & 0 &$\Bbbb_{32}^T$& $\Bbbb_{42}^T$& 0\\\hline\hline
$3(k+2)n_\cE$ & 0 &$\Bbbb_{32}$& 0 & 0 & $\Fb_{3}$\\\hline \hline
$3(k+4)n_\cE+6n_\cV$ &$\Bbbb_{41}$&$\Bbbb_{42}$ & 0 & 0 & 0\\\hline\hline
%&$3(k+1)n_\cE$& $3(k+1)n_\cE$& $3n_\cE$& $3n_\cE$ & $3n_\cE$ & $3n_\cE$ & $3$ & $3$ & $3(k+2)n_\cE$ & $3(k+2)n_\cE$ & $3n_\cV$ & $3n_\cV$ &
\end{tabular}
\end{center}
or in more detail, %ed structure of the matrix $\Kbb$.
\begin{center}
%\begin{tabular}{c||c|c|c|c|c|c|c|c||c||c|c|c||c}
%&$\qb$ & $\pb$ & $\Pb_+$ & $\Pb_-$ & $\Qb_+$ & $\Qb_-$ & $\alb$ & $\betab$ & $\ub$ & $\omegab$ & $\Ub$ & $\Omegab$& dimension\\\hline\hline
%$\xib$ &\xx & &&&&&&&0& \xx & && $3(k+1)n_\cE$\\\hline \hline
%$\thb$ & &&&&&&&& \xx & \xx & && $3(k+1)n_\cE$\\\hline
%$\Thb_+$ &&&&&&&&&\xx &&$-(\Abb_\cI^+)^T$&& $3n_\cE$\\\hline
%$\Thb_-$ &&&&&&&&&\xx &&$(\Abb_\cI^-)^T$&& $3n_\cE$\\\hline
%$\Xib_+$ &&&&&&&&&&\xx &&$-(\Abb_\cI^+)^T$& $3n_\cE$\\\hline
%$\Xib_-$ &&&&&&&&&&\xx &&$(\Abb_\cI^-)^T$& $3n_\cE$\\\hline
%$\gamb$ &&&&&&&&&\xx &&&& $3$\\\hline
%$\deltab$ &&&&&&&&&\xx &&&& $3$\\\hline\hline
%$\vb$ &0&\xx&\xx&\xx&&&\xx&&0&&&& $3(k+2)n_\cE$\\\hline \hline
%$\omegab$ &\xx&\xx&&&\xx&\xx&&\xx&&&&& $3(k+2)n_\cE$\\\hline
%$\Vb$ &&&$-\Abb_\cI^+$&$\Abb_\cI^-$&&&&&&&&& $3n_\cV$\\\hline
%$\Wb$ &&&&&$-\Abb_\cI^+$&$\Abb_\cI^-$&&&&&&& $3n_\cV$\\\hline\hline
%%&$3(k+1)n_\cE$& $3(k+1)n_\cE$& $3n_\cE$& $3n_\cE$ & $3n_\cE$ & $3n_\cE$ & $3$ & $3$ & $3(k+2)n_\cE$ & $3(k+2)n_\cE$ & $3n_\cV$ & $3n_\cV$ &
%\end{tabular}
\includegraphics[bb=70 565 540 760,width=\textwidth]{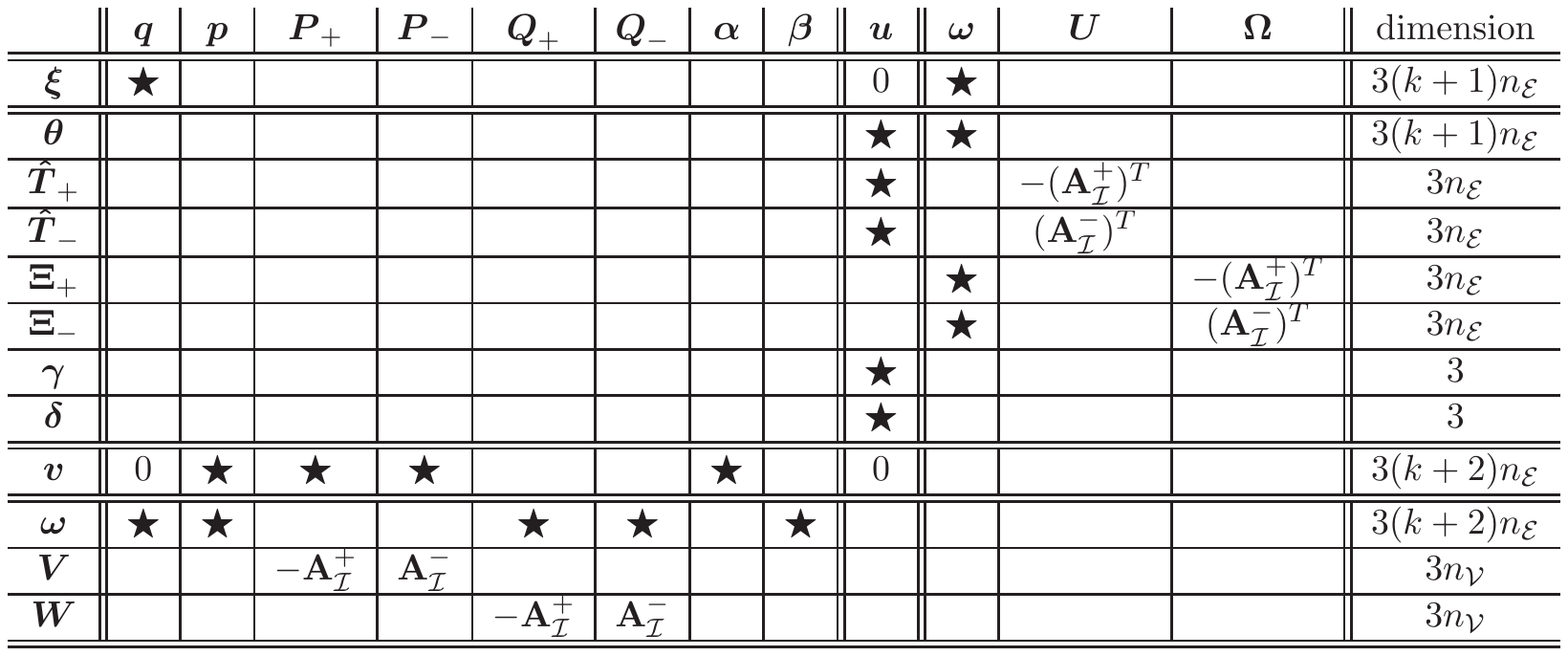}
\end{center}

\subsection{Numerical results}

\begin{figure}[!]
\begin{center}
\includegraphics[bb=158 308 510 479,width=0.4\textwidth]{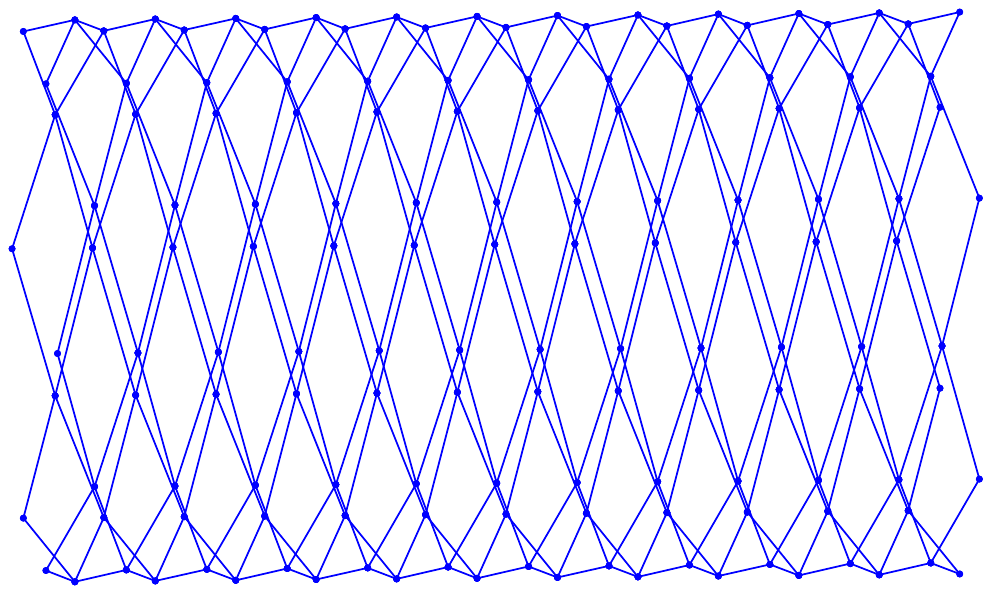}\quad
\includegraphics[bb=179 301 421 496,width=0.3\textwidth]{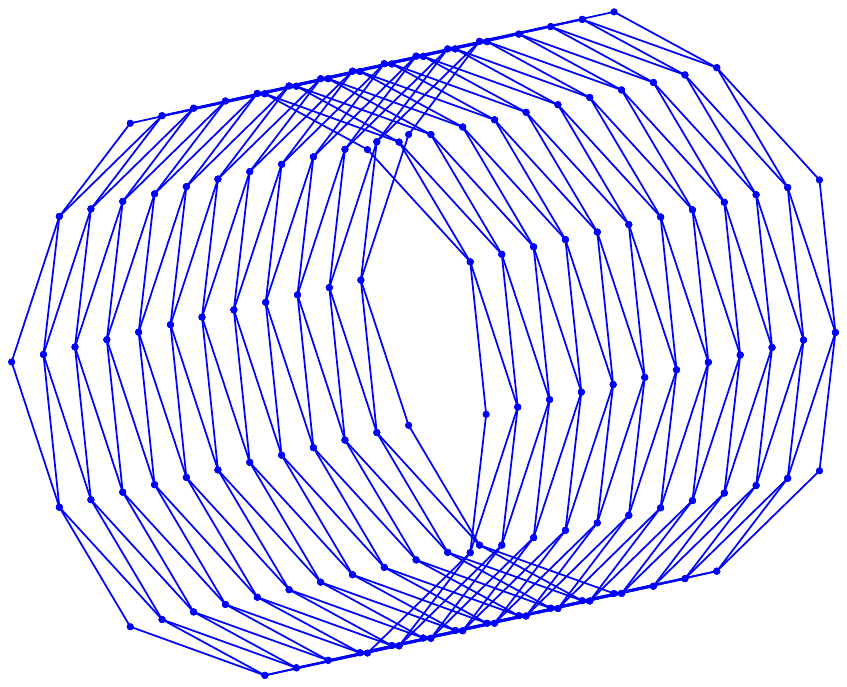}
\caption{Design of the Palmaz like stent used in simulations\label{fPalmaz}}
\end{center}
\end{figure}

To illustrate our theoretical analysis we test the implementation of the numerical scheme in the new formulation for a Palmaz type stent as in Figure~\ref{fPalmaz}.
The radius of the stent is $1.5$mm and the overall length is  $1.68$cm. There are $144$ vertices in the associated graph with $276$ straight edges. All vertices except the boundary ones are junctions of four edges.
The cross-sections are assumed to be square with the side length $0.1$mm. The material of the stent is stainless steel with Young modulus $E=2.1 \cdot 10^{11}$ and Poisson ratio $\nu=0.26506$.
To this structure we apply the forcing normal to the axis of the of stent, i.e. of the form
\begin{equation}\label{radialf}
\fb (\xb) =  f(x_1) \frac{x_2\eb_2+x_3\eb_3}{\sqrt{x_2^2+x_3^2}}, \qquad \xb =(x_1,x_2,x_3) \in \ZR^3,
\end{equation}
where $x_1$ is the axis of the cylinder. As a consequence, the deformation will also posess some radial symmetry.
The problem is a pure traction problem and the applied forces satisfy the necessary condition. The non-uniqueness of the solution in the problem is fixed using the Lagrange multipliers $\alb$ and $\betab$.

\begin{figure}[h!]
\begin{center}
\includegraphics[bb=96   238   515   553,width=0.44\textwidth]{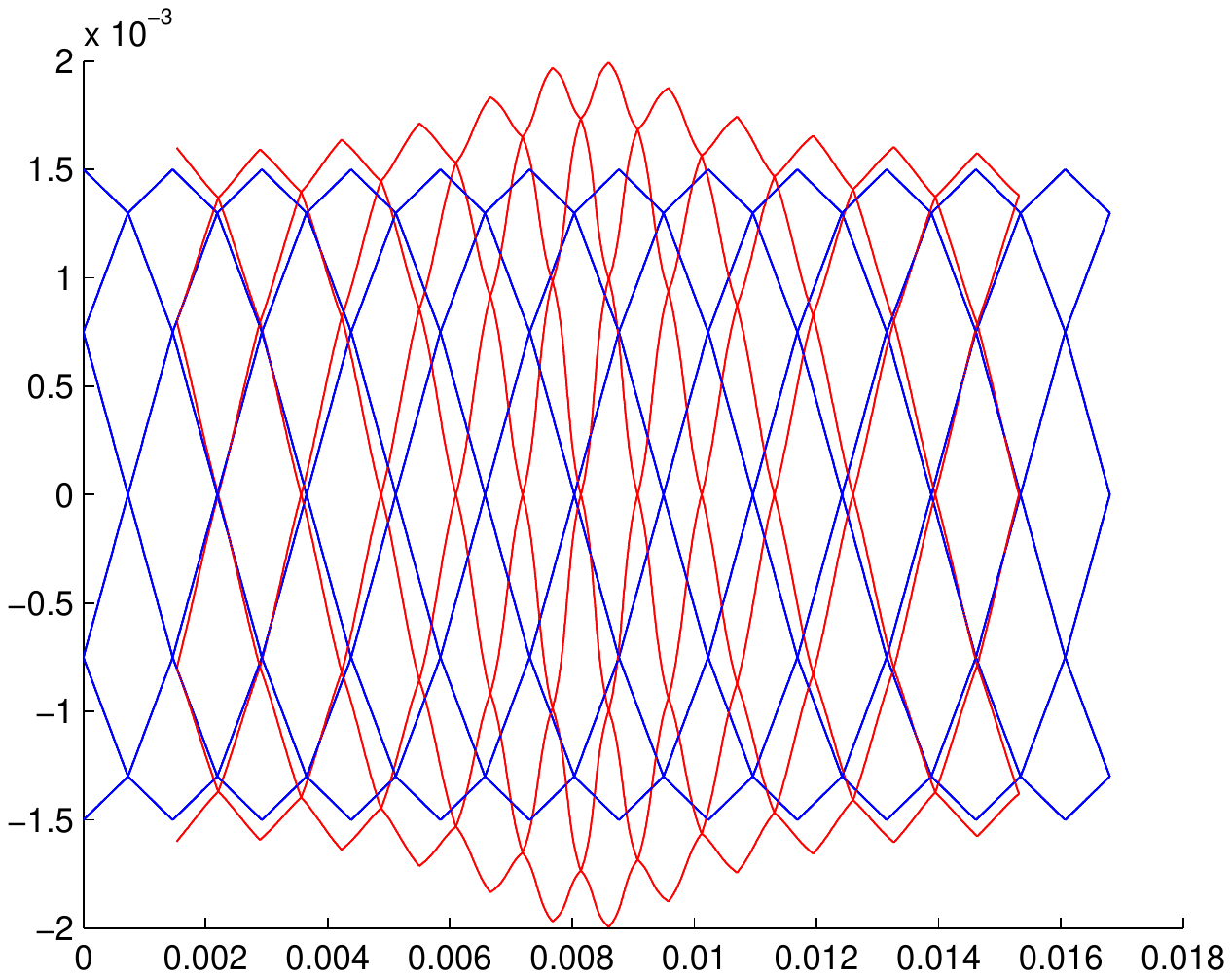}\quad
\includegraphics[bb=96   238   515   553,width=0.5\textwidth]{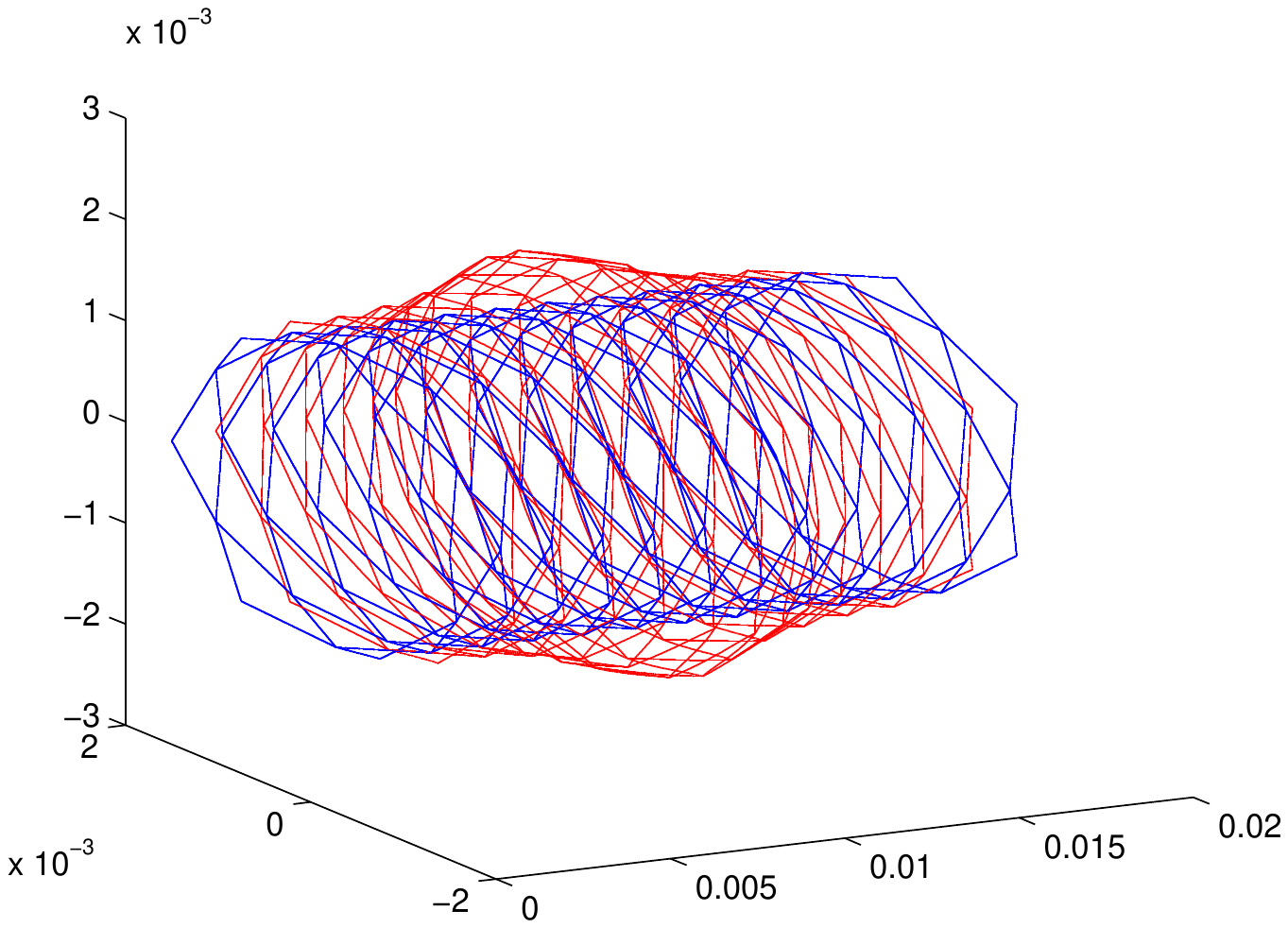}
\caption{Solution for the radial load from \eqref{f1}. On the left the solution is projected to the $(x_1,x_2)$--plane; on the  right the solution is shown from a different perspective.}\label{fPalmazrj}
\end{center}
\end{figure}

\begin{figure}[h!]
\begin{center}
\includegraphics[bb=96   238   515   553,width=0.5\textwidth]{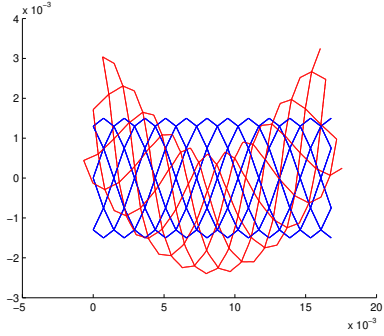}
\caption{Solution for the load given in \eqref{f2} projected into $(x_1,x_2)$--plane}\label{bended_stent}
\end{center}
\end{figure}

The solution for the forcing function
\begin{equation}\label{f1}
f(x_1) = \frac{10}{10^5(x_1-\ell/2)^2+1},
\end{equation}
is presented in Figure~\ref{fPalmazrj}; here $\ell$ is the length of the stent. On the left the solution is projected to the $(x_1,x_2)$--plane, while on the  right it is shown from the different perspective.
For the forcing function
\begin{equation}\label{f2}
\fb(x_1) = 10^3 (x_1-\ell/2)^2\eb_3
\end{equation}
the results are given in Figure~\ref{bended_stent}; again $\ell$ is the length of the stent.

%
%while for quadsatic
%$$
%f(x_1) = 25000000 x_1^2 NOVA FUNKCIJA!
%$$
%the solution is plotted in Figure~\ref{fPalmazrj2}.
%\begin{figure}[h!]
%\begin{center}
%\includegraphics[width=0.4\textwidth]{rj1x2.eps}\quad
%\includegraphics[width=0.45\textwidth]{rj2x2.eps}
%\caption{Solution for the radial force quadsatic with respect to the length of the stent}\label{fPalmazrj2}
%\end{center}
%\end{figure}
In the following we present the order of convergence of the finite element method for the solution of the problem with the quadratic forcing $f(x_1) = 2.5 \cdot 10^7 x_1^2$, the same as in the numerical scheme presented in \cite{GT}.
We divide all the edges into $128$ smaller rods and solve the equilibrium problem. The obtained solution we consider as the best possible and use it to compute the errors, denoted by "$\mbox{error}(i)$", of the approximations for edges split into  $2^i$ smaller struts, $i=1,\ldots,6$. We use {quadratic finite elements for displacement and infinitesimal rotation, and linear finite elements for contact forces and couples, i.e. $n=2$ and $k=1$}, see e.g. \cite{CiarletFEM}, for computing the approximations, the $L^2$ norm and the $H^1$ semi-norm for displacements and the $\ell^1/n$ norm for unknowns in $\ZR^n$, i.e., the arithmetic mean of errors, to determine the error estimates and these to
compute the convergence rates via
\begin{equation}\label{cr7}
\frac{\log{\frac{\mbox{error}(i+1)}{\mbox{error}(i)}}}{\log{\frac{h(i+1)}{h(i)}}}, \qquad i = 1, \ldots, 5.
\end{equation}
\begin{figure}[h!]
\begin{center}
\includegraphics[bb=96   238   515   553,width=0.6\textwidth]{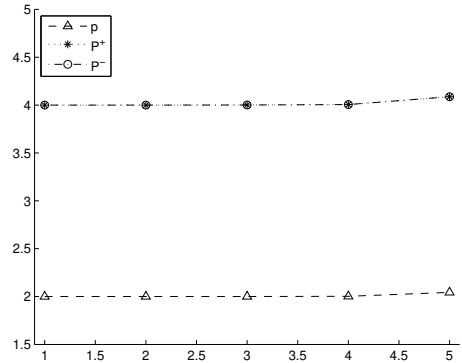}
\caption{Rate of convergence (\ref{cr7}) for the $L^2$ norm of the contact force $\pb$ and $\ell^1$ norm of $\Pb_+,\Pb_-$\label{redkvgP}}
\end{center}
\end{figure}
\begin{figure}[h!]
\begin{center}
\includegraphics[bb=96   238   515   553,width=0.6\textwidth]{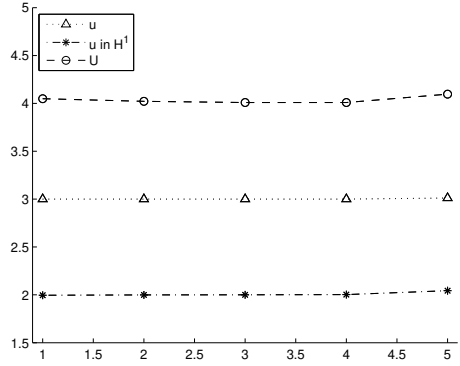}
\caption{Rate of convergence (\ref{cr7}) for the $H^1$ seminorm and $L^2$ norm of the displacement $\ub$ and $\ell^1$ norm of $\Ub$\label{redkvgU}}
\end{center}
\end{figure}
The obtained convergence rates for $\ub,\omegab,\pb,\qb$ are in agreement with the analytical estimate from Theorem~\ref{terrorglavni} for $k=1$ and $n=2$. {In Figure~\ref{redkvgP} the convergence rates for $\pb,\Pb_+$ and $\Pb_-$ are displayed, while in Figure~\ref{redkvgU} the convergence rates for the $\ub$ ($L^2$ norm and $H^1$ semi-norm) and $\Ub$ are plotted.
Additionally we present the errors of the remaining unknowns in Table~\ref{table2}.
\begin{table}[h!]
\begin{center}
\begin{tabular}{r|rrrrr}
splitting&$q$ & $Q_+$ & $Q_-$ & $\omega$ & $\Omega$\\\hline
 2&1.8127e-5 & 3.1945266e-8 & 3.2871209e-8 & 8.52118e-4 & 3.3700235e-5\\
 4&0.4532e-5 & 0.1995776e-8 & 0.2024711e-8 & 1.06516e-4 & 0.2116646e-5\\
 8&0.1133e-5 & 0.0124898e-8 & 0.0125802e-8 & 0.13314e-4 & 0.0132275e-5\\
16&0.0283e-5 & 0.0007811e-8 & 0.0007839e-8 & 0.01664e-4 & 0.0008260e-5\\
32&0.0070e-5 & 0.0000486e-8 & 0.0000487e-8 & 0.00208e-4 & 0.0000514e-5\\
64&0.0017e-5 & 0.0000028e-8 & 0.0000028e-8 & 0.00025e-4 & 0.0000030e-5
\end{tabular}
\caption{Errors for different splitting of edges ($L^2$ errors for $q$ and $\omega$).\label{table2}}
\end{center}
\end{table}

To compare the numerical scheme for the new formulation with that for the old formulation in \cite{GT} in Table~\ref{table1}, we present the obtained matrix sizes and the computing times (computing times are determined using difference of "toc" and "tic" functions in MATLAB 2010b.)  for the computation on a personal computer with 16GB RAM, 64bit Windows and with INTEL\textregistered Core i3-7100 CPU@3.90FHz.
\begin{table}[h!]
\begin{center}
\begin{tabular}{c|c|c|c|c}
& \multicolumn{2}{c|}{new formulation}	&	\multicolumn{2}{c|}{old formulation}	\\ \hline
splitting $\#$ 	&  comp. time in $s$	 &size of matrix	&comp. time in $s$	&size of matrix\\\hline
8		&22	&105198	&2	&38958\\
16		&47	&211182	&42	&78702\\
32		&108	&423150	&152	&158190\\
64		&288	&847086	&629	&317166\\
128		&903	&1694958	&4183	&635118\\
\end{tabular}
\caption{Computing times and matrix sizes for the old and new numerical scheme\label{table1}}
\end{center}
\end{table}

The difference between the solution for the displacement $\ub$, infinitesimal rotation $\omegab$ and the contact force $\pb$ (the only quantities calculated in the old numerical scheme from \cite{GT}) for the same splitting of edges is
at least four digits smaller than the error of the approximation, i.e. we obtain the same order of approximation for the same mesh. Further, the size of the matrix for the new numerics is 2.7 times larger. However, the computing times for the new approach are smaller when splitting the edges in $32$ or more rods.

\section{Dynamic modeling  of elastic stents}
\setcounter{equation}{0}
\label{evol}

In the previous sections we have considered the static problem for the stent, but for the analysis, in particular, to study the movement of the stent under the permanent excitation through the heartbeat, in this section we formulate and analyze the evolution model of the stents.
% and the associated eigenvalue problem.

\subsection{Formulation of the space-continuous dynamic model}

In order to formulate the dynamic model we start from the evolution equation of curved rods from \cite{Tevolution} and add in  the model (\ref{lj1ri})--(\ref{lj4ri}) the inertial term to the equilibrium equation (\ref{lj1ri}), which changes to
\[
\rho A \ub^i_{tt} = \partial_s \pb^i + \fb^i,
\]
where $A$ is the area of the cross section and $\rho$ is volume density of mass.
In the weak formulation this implies that the term
\[
-\sum_{i=1}^{n_\cE}\int_0^{\ell^i} \rho A \ub^i_{tt} \cdot \vb^i\, ds
\]
should be added to the left hand side of (\ref{lj12}) and (\ref{lj12A}). Using the notation of Section~\ref{newmodel} and introducing the bilinear form
\[
c: M \times M \to \ZR, \qquad c(\phib,\psib) = \sum_{i=1}^{n_\cE}\int_0^{\ell^i} \rho A \ub^i \cdot \vb^i\, ds,
\]
we can formulate the evolution problem of elastic stents as follows:

Determine $\funS$ and $\funG$ such that
\begin{equation}\label{evoeqmain}
\aligned
a(\funS,\funG) + b(\funG,\funu) &= 0, \qquad\qquad \funG \in V,\\
-\frac{d^2}{dt^2}c(\funu,\funv) + b(\funS,\funv) &= f(\funv), \qquad \funv \in M.
\endaligned
\end{equation}
%

%
%\subsection{Analysis of the continuous model}
%%
%For this we will first analyze the spectral properties of the operator pencil associated to the system \eqref{evoeqmain}, see \cite{xxx} for the general theory.
%\marginpar{Menicken/M\"{o}ller or another reference}

\subsection{Analysis of the space-discrete model}
The dynamical system \eqref{disdynsys} is a finite dimensional  linear differential-algebraic equation of second order, with a nonsingular (but indefinite) stiffness matrix $\Kbb$, and a positive semidefinite but singular $\Ebb$.

The associated discrete dynamical problem is given by
\begin{equation} \label{disdynsys}
-\Ebb \ddot{\zb}(t) + \Kbb \zb(t) = \Fb(t),
\end{equation}
where the matrix $\Kbb$ and the right hand side $\Fb$ are as in the static case, $\Ebb$ has the following structure partitioned as $\Kbb$
\[
%\begin{center}
\begin{tabular}{c||c||c||c||c||}
dim $\backslash$ unknown&$\qb$ & $(\pb,\Pb_+,\Pb_-,\Qb_+,\Qb_-,\alb,\betab)$ & $\ub$ & $(\omegab,\Ub,\Omegab)$ \\\hline\hline
$3(k+1)n_\cE$ &$0$ & 0 & 0 & 0\\\hline \hline
$3(k+5)n_\cE+6$ & 0 & 0 &0& 0\\\hline\hline
$3(k+2)n_\cE$ & 0 &0& \Mbb & 0 \\\hline \hline
$3(k+4)n_\cE+6n_\cV$ &0&0 & 0 & 0 \\\hline\hline
%&$3(k+1)n_\cE$& $3(k+1)n_\cE$& $3n_\cE$& $3n_\cE$ & $3n_\cE$ & $3n_\cE$ & $3$ & $3$ & $3(k+2)n_\cE$ & $3(k+2)n_\cE$ & $3n_\cV$ & $3n_\cV$ &
\end{tabular},
\]
%\end{center}
and $\zb(t)$ is the vector of coefficient functions in the finite element basis.}

The particular block structure of the DAE allows to analyze the properties of the system, which are characterized by the spectral properties of the matrix polynomial $-\lambda^2 \Ebb +\Kbb$.
\begin{lemma}\label{canform}
Consider the  DAE \eqref{disdynsys} and the associated pair of matrices $(\Ebb,\Kbb)$. Then there exists a nonsingular matrix $\Vbb$ with the property that
\begin{eqnarray*}
\hat \Ebb&=&\Vbb^T \Ebb \Vbb=\mat{ccccc} 0 &&&& \\ & 0 &&& \\ && 0 &&\\ &&& \Mbb & \\ &&&&0 \rix,\quad \Vbb^{-1} \zb= \mat{c} \hat \zb_1 \\ \hat \zb_2\\ \hat \zb_3 \\ \hat \zb_4 \\ \hat \zb_5 \rix, \\
\hat \Kbb&=&\Vbb^T \Kbb \Vbb= \mat{ccccc} 0 & 0 & 0 & 0 & \hat \Bbbb_{51}^T \\
0 & \hat \Abb_{22} & 0 & \hat \Bbbb_{42}^T & 0 \\
0 & 0 & \hat \Abb_{33} & 0 & 0\\
0 & \hat \Bbbb_{42} & 0 & 0 & 0 \\
\hat \Bbbb_{51} & 0 &  0 & 0 & 0 \rix,\quad \quad \Vbb^T \Fb= \mat{c} 0\\ 0\\ 0 \\ \hat \Fb_4\\ 0\rix,
\end{eqnarray*}
where $\hat \Abb_{33}=\hat \Abb_{33}^T$, $\hat \Bbbb_{42}$, and $\hat \Bbbb_{51}$ are invertible, and $\hat \Fb_4=\Fb_3$.
\end{lemma}
\proof
The proof follows by a sequence of congruence transformations, starting from  the original block structure
\[
\Ebb=\mat{ccccc} 0 &&& \\ & 0 && \\  && \Mbb& \\ &&&0 \rix,\quad  \Kbb= \mat{ccccc} \Abb_{11} & 0 &  0 & \Bbbb_{41}^T \\
0 & 0 & \Bbbb_{32}^T & \Bbbb_{42}^T \\
0 & \Bbbb_{32} & 0 & 0 \\
\Bbbb_{41} & \Bbbb_{42} & 0 & 0 \rix,
\]
by first compressing
\[\mat{cc} 0 & \Bbbb_{32}  \\
\Bbbb_{41} & \Bbbb_{42} \rix
\]
via an orthogonal transformation $\Vbb_1$ from the right to a form
\[
\mat{ccc}
\tilde \Bbbb_{41} & \tilde \Bbbb_{42} & 0 \\
\tilde \Bbbb_{51} & 0 &  0  \rix,
\]
row-partitioned according to the original row-partitioning,
with $\tilde \Bbbb_{42}$, and $\tilde \Bbbb_{51}$ having full column rank.  Setting
$\tilde \Vbb=\diag( \Vbb_1, \Ibb)$ and applying the transformation with $\Vbb_1^T$ from the right to the first two block columns, with $\Vbb_1$ from the left to the first two block rows, and partitioning
\[
\Vbb_1^T \mat{cc} \Abb_{11} & 0 \\ 0 & 0 \rix\Vbb_1=: \mat{ccc} \tilde \Abb_{11}  &\tilde \Abb_{12} &\tilde \Abb_{13}\\ \tilde \Abb_{21} &\tilde \Abb_{22} &\tilde \Abb_{23} \\\tilde \Abb_{31} & \tilde \Abb_{32} &\tilde \Abb_{33} \rix
\]
accordingly, we obtain a transformed system with
\begin{eqnarray*}
\tilde \Ebb&=&\tilde \Vbb^T \Ebb \tilde \Vbb =\mat{ccccc} 0 &&&& \\ & 0 &&& \\ && 0 &&\\ &&& \Mbb & \\ &&&&0 \rix,\quad {\tilde \Vbb}^{-1} \zb= \mat{c} \tilde \zb_1 \\ \tilde \zb_2\\ \tilde \zb_3 \\ \tilde \zb_4 \\ \tilde \zb_5 \rix, \\
\tilde \Kbb&=&\tilde \Vbb^T \Kbb \tilde \Vbb= \mat{ccccc}  \tilde \Abb_{11}  &\tilde \Abb_{12} &\tilde \Abb_{13}& \tilde \Bbbb_{41}^T & \tilde \Bbbb_{51}^T \\ \tilde \Abb_{21} &\tilde \Abb_{22} &\tilde \Abb_{23} & \tilde \Bbbb_{42}^T & 0 \\
\tilde \Abb_{31} & \tilde \Abb_{32} &\tilde \Abb_{33}  & 0 & 0\\
\tilde \Bbbb_{41} & \tilde \Bbbb_{42} & 0 & 0 & 0 \\
\tilde \Bbbb_{51} & 0 &  0 & 0 & 0 \rix,\quad \tilde \Fb:=\tilde \Vbb^T \Fb= \mat{c} 0\\ 0\\ 0 \\ \hat \Fb_4\\ 0\rix.
\end{eqnarray*}
The property that $\Abb$ is $\Ker \Bbbb$ elliptic then implies that $\tilde \Abb_{33}$ is nonsingular. The fact that
$\Kbb$ is invertible and that $\tilde \Bbbb_{51}$ has full column rank implies that $\tilde \Bbbb_{51}$ is square and nonsingular. Thus, there exists a nonsingular matrix $\Vbb_2$ that eliminates the leading $4$ blocks in the first block row and column with $\tilde \Bbbb_{51}$ and the leading $2$ blocks in the third block column and column with $\tilde \Abb_{33}$. Thus, a congruence transformation with $\Vbb= \Vbb_2 \tilde \Vbb$ yields the asserted structure. The property that $\tilde \Bbbb_{42}$ has full column rank and that $\hat\Kbb$ is invertible, then implies that $\hat \Bbbb_{42}=\tilde \Bbbb_{42}$ is square and nonsingular as well.
\eproof

\begin{remark}{\rm
The elimination procedure presented in the proof of Lemma~\ref{canform} is a structured version of the canonical form construction for differential-algebraic equations, see e.g. \cite{KunM06}, which identifies all explicit or implicit constraints in the system. These lead to restrictions in the initial values and usually also to further differentiability conditions for the inhomogeneity. Due to the structure of the equations and the inhomogeneity, the only components of the inhomogeneity that have to be differentiated in time, are $0$, so
there are no further requirements on the forcing function $\Fb$.

Note further that the same procedure can also be applied in the space-continuous case, by constructing appropriate projections into subspaces, see e.g. \cite{LamMT13}. Since this procedure is rather technical we do not present this construction here.}
\end{remark}

\begin{corollary}\label{solution}
Consider the  DAE \eqref{disdynsys} transformed as in Lemma~\ref{canform}.
Then for the general solution of the transformed system we have $\hat \zb_1=0$, $\hat \zb_3=0$, $\hat \zb_5=0$ ,
$\hat \zb_2$ is the solution of the second order DAE
\begin{equation}\label{z2DAE}
\hat \Abb_{22}\ddot {\hat \zb}_2 =  \hat \Bbbb_{42}^T \Mbb^{-1} \hat \Bbbb_{42} {\hat \zb}_2  + \hat \Bbbb_{42}^T \Mbb^{-1}  \hat \Fb_4
\end{equation}
which exists without any further smoothness requirements for $\hat \Fb_4$, and is unique for every consistent initial  condition, and finally $\hat \zb_4=  -\hat \Bbbb_{42}^{-T} \hat \Abb_{22} {\hat \zb}_2$.

No initial conditions can be assigned for $\hat \zb_1$, $\hat \zb_3$, $\hat \zb_5$, $\hat \zb_4$, while for $\hat \zb_2$ a consistency condition for $\dot{\hat \zb}_2$ in the kernel of $\hat \Abb_{22}$ arises, that depends on the right hand side $\hat \Fb_4$.
\end{corollary}
\proof
This follows directly from the transformed equations. The equations (\ref{z2DAE}) form a so called index-one DAE (see \cite{KunM06}), since  $\hat \Bbbb_{42}^T \Mbb^{-1} \hat \Bbbb_{42}$ is positive definite and thus, in particular invertible  in the kernel of $\hat \Abb_{22}$. Projecting  into this kernel gives an algebraic equation which has to hold for the initial condition associated with $\dot {\hat \zb}_2$, while for the remaining components of $\hat \zb_2$ an initial value can can be chosen arbitrarily.
\eproof

\begin{remark}{\rm
	The operator pencil associated to the system \eqref{evoeqmain} can be studied using the general theory from \cite{MenM03}. Note that the form $e( (\funS,\funu), (\funG,\funv))=c(\funu,\funv)$, $(\funS,\funu),(\funG,\funv)\in V\otimes M$ is bounded and symmetric and the form $k((\funS,\funu), (\funG,\funv)) = a(\funS,\funG) + b(\funG,\funu) + b(\funS,\funv)$, $(\funS,\funu),(\funG,\funv)\in V\otimes M$ is closed, symmetric and semi-bounded from below, see \cite{GKMV}. The form $e$ defines, in the sense of Kato \cite{Kat13}, the bounded, semidefinite and self-adjoint operator $E$, whereas the form $k$ defines the self-adjoint semibounded from below operator $K$. The operator $K$ can be represented as the formal product $K=L^*JL$, where $L$ is a closed operator with a bounded inverse such that the domains of $L$ and $k$ are equal and $J$ is the so called fundamental symmetry (a bounded self adjoint operator such that $J^2=J$). Subsequently it can be shown --- using the special structure of $E$ --- that the operator function $T(z)=-z^2\tilde{E}+\tilde{K}$, where $\tilde{E}=L^{-*}EL^{-1}$ and $\tilde{K}=J$, is Fredholm operator valued. Furthermore, the operators $\tilde{E}$ and $\tilde{K}$ are bounded Hermitian operators and the resolvent set of $T$ contains zero by Theorem \ref{main}. By the results of \cite[Section 1]{MenM03}, the pencil $T$ has finite semi-simple eigenvalues of finite multiplicity. The construction of an oblique projection onto the reducing subspace associated to the eigenvalue infinity can be done in a similar way as in Lemma \ref{canform}. The construction is decidedly more technical and we leave it to a subsequent paper where we will discuss more general second order systems (e.g. those involving a damping term).}
\end{remark}

\subsection{Numerical results}

In this subsection we present some numerical results obtained for the evolution problem.
The time discretization is done using the implicit mid point rule \cite{HaiW96} (of convergence order $2$)  applied to the first order formulation of the system.
At each time step a linear system for the matrix $-\Ebb+0.25  \Delta t^2 \Kbb$ is solved using backslash in MATLAB. The computations are performed on a personal computer with 16GB RAM, 64bit Windows and with INTEL\textregistered Core i3-7100 CPU@3.90FHz. The presented computing times are determined using the difference of "toc" and "tic" functions in MATLAB 2010b.
All simulations are carried out for the following set of parameters:
\begin{itemize}
\item elasticity coefficients: $\mu = E = 1$ Pa,
\item thickness of the stent struts: $0.0001$ m,
\item load: radial, as given in (\ref{radialf}), where
\begin{equation}\label{force}
f(x,t) = F\left(\frac{\pi}{0.003}\left(x-c_{\textrm{vawe}}(t-t_0) \right)\right)
\end{equation}
with
\[
F(y) = \left\{\begin{array}{ll}
5 \cdot 10^{-8} \cos(y), &  \textrm{ if } |y|< 0.0015\\
0, & \textrm{ else}
\end{array}
\right.,  \ c_{\textrm{vawe}}=0.0075, \ t_0 = 0.5.
\]
In other words, $f$ is given as traveling wave determined by the function $F$, where the factor $c_{\textrm{vawe}}$ denotes the speed of the wave and the term $t_0$ asserts the condition $f(x,\cdot)=0$.
\item mass density $\rho = 2000\; kg/m^3$,
\item total time $T=12$s.
\end{itemize}

\begin{figure}[h!]
\begin{center}
\includegraphics[bb = 96   238   515   553, width=0.45\textwidth]{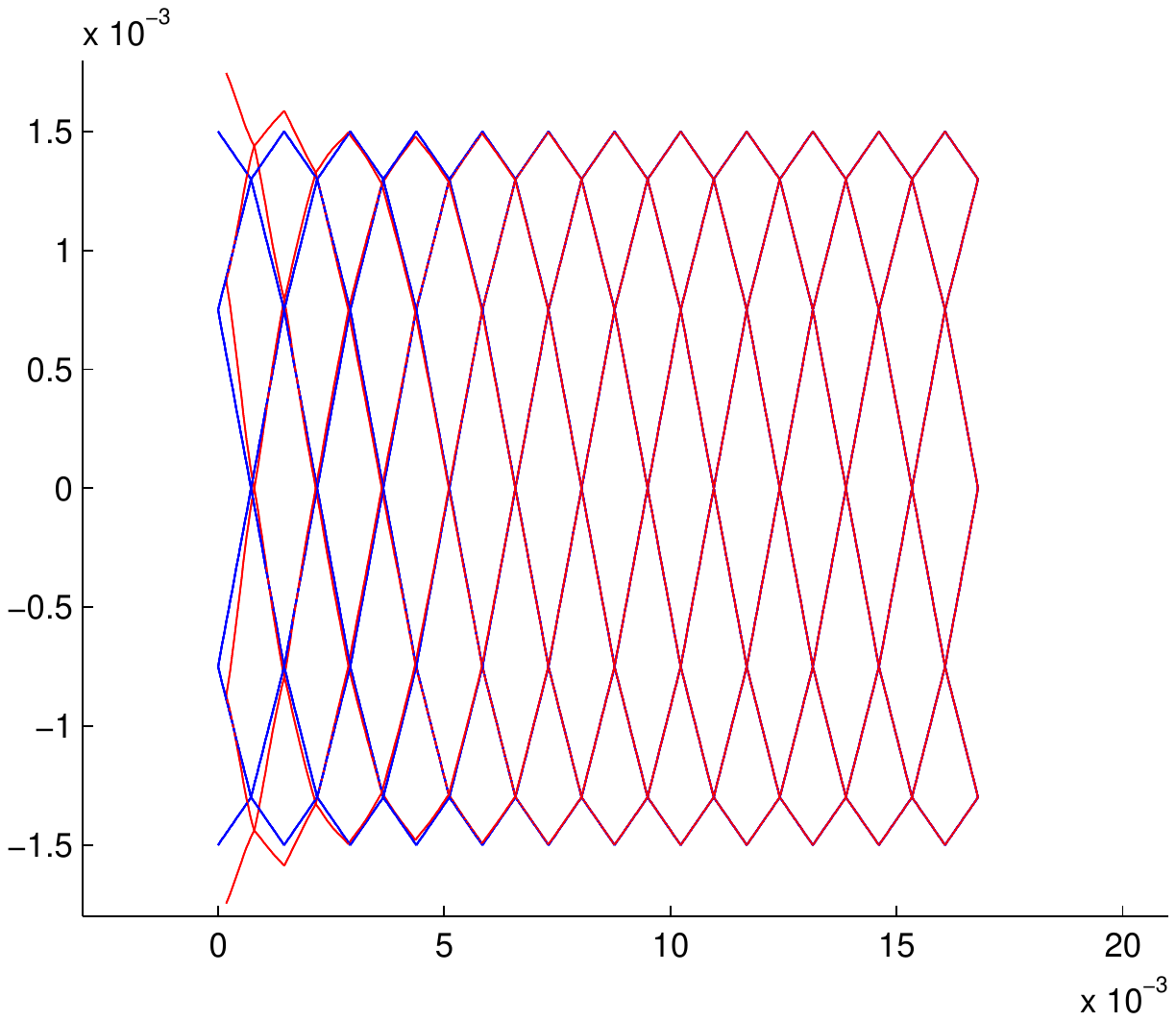}
\qquad
\includegraphics[bb = 96   238   515   553, width=0.45\textwidth]{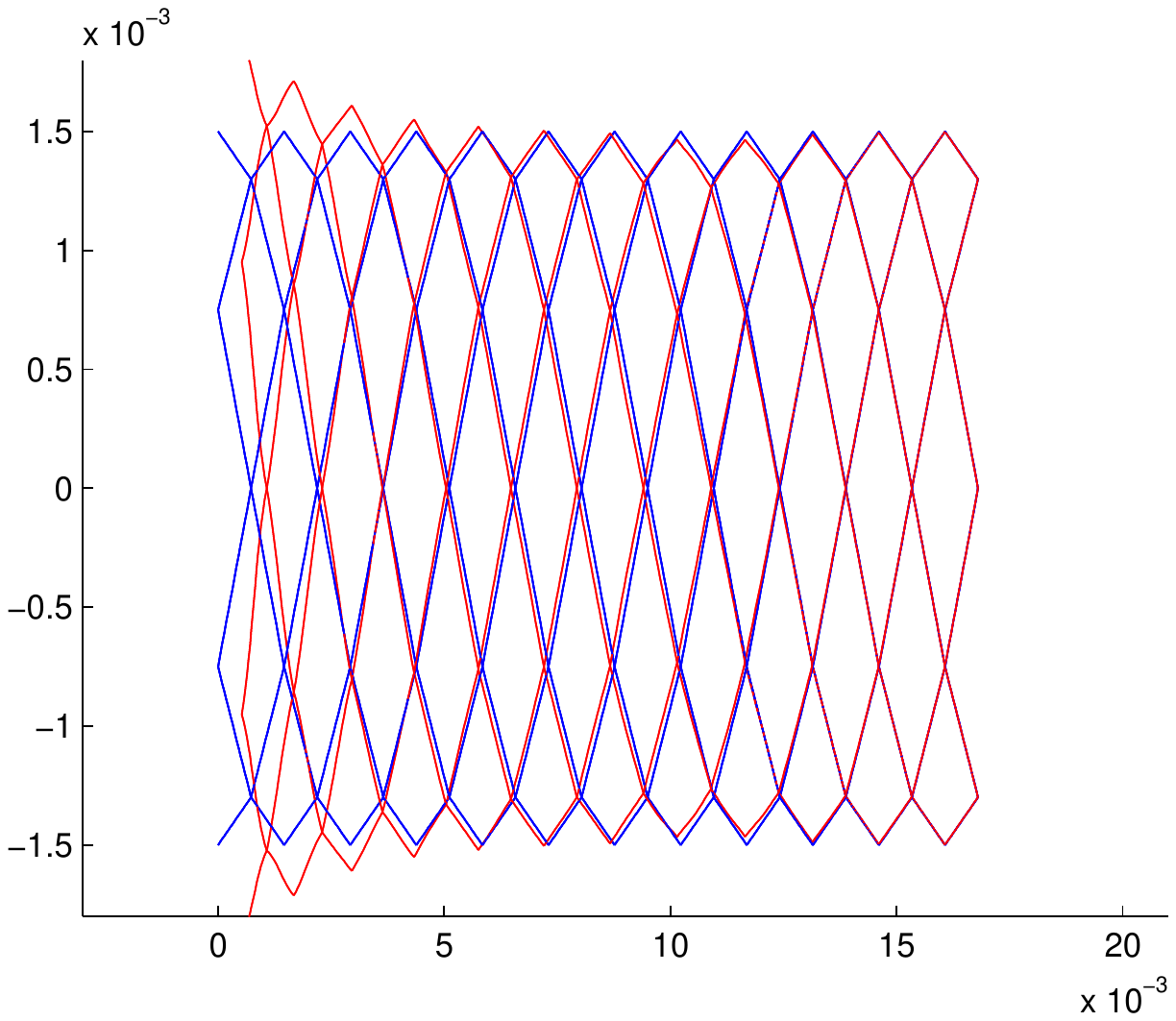}\\
\includegraphics[bb = 96   238   515   553, width=0.45\textwidth]{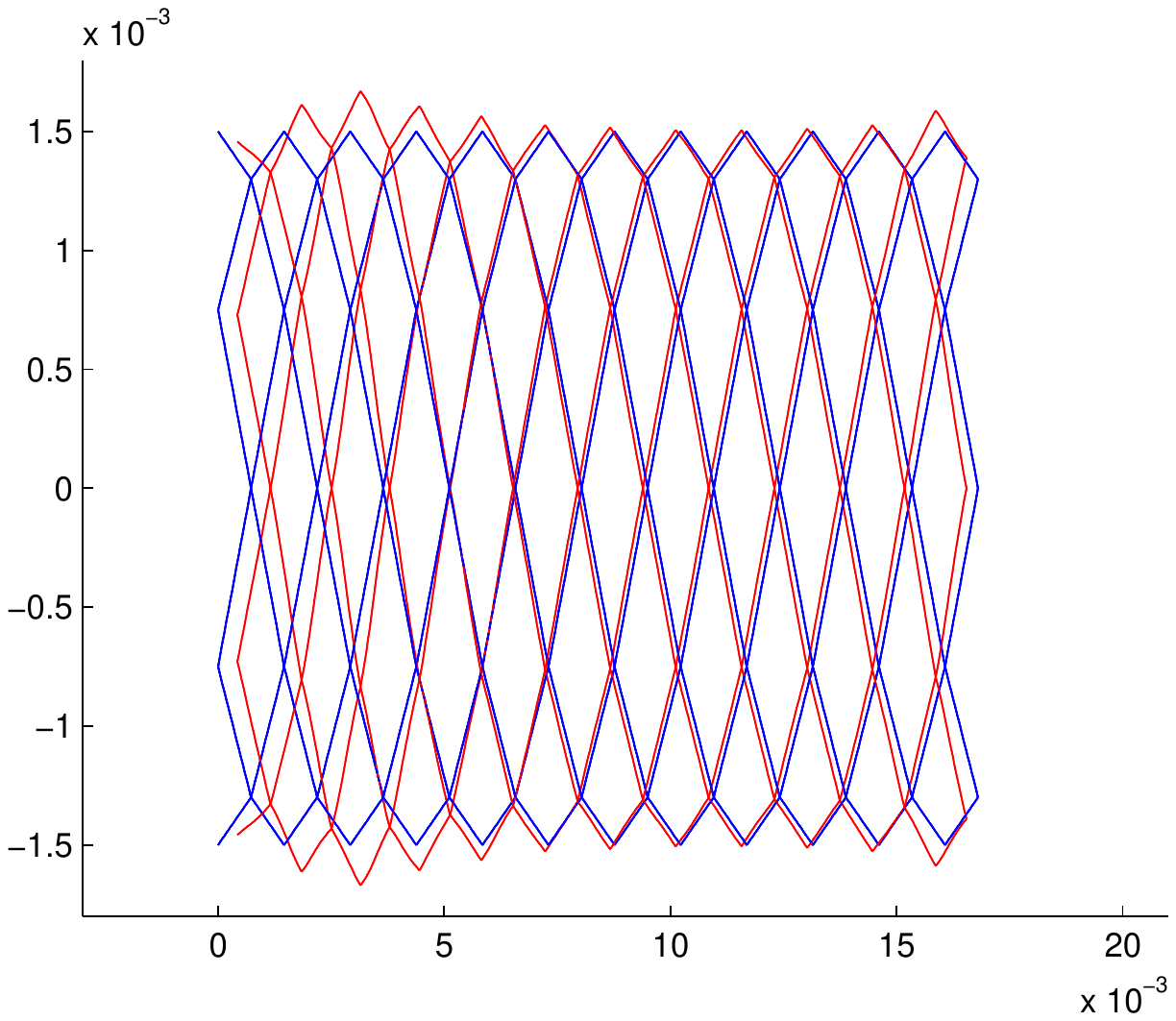}
\qquad
\includegraphics[bb = 96   238   515   553, width=0.45\textwidth]{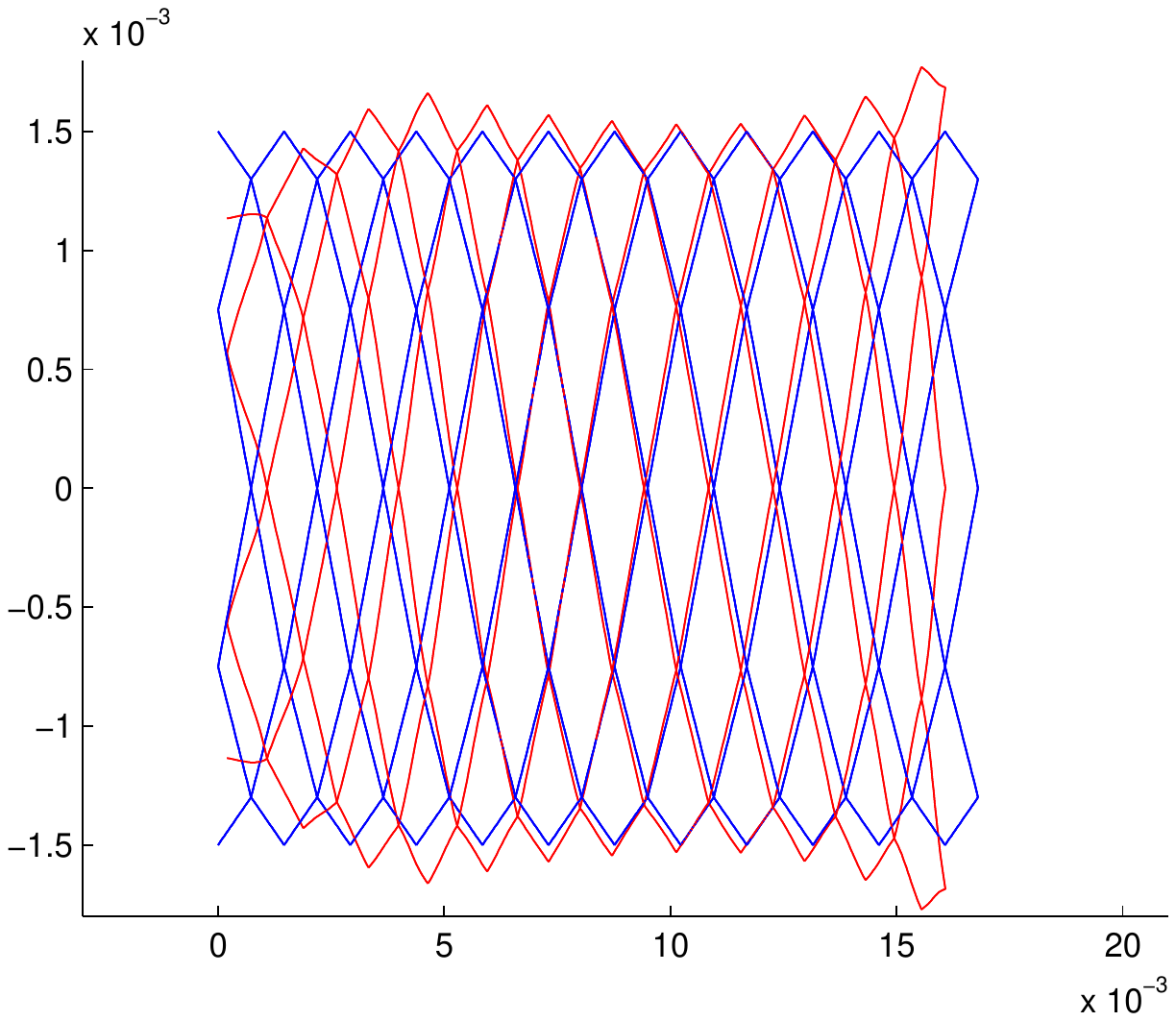}\\
\includegraphics[bb = 96   238   515   553, width=0.45\textwidth]{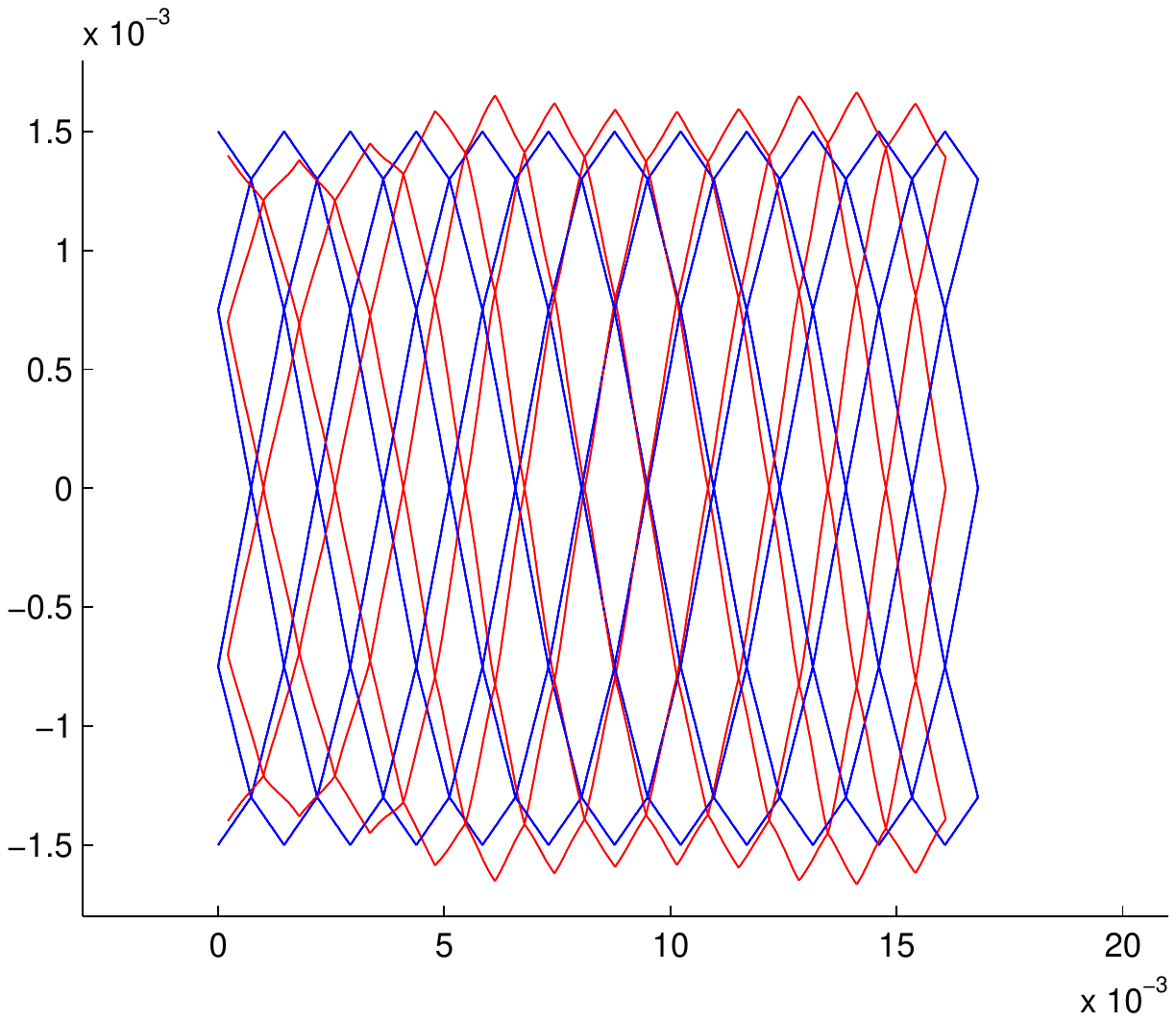}
\qquad
\includegraphics[bb = 96   238   515   553, width=0.45\textwidth]{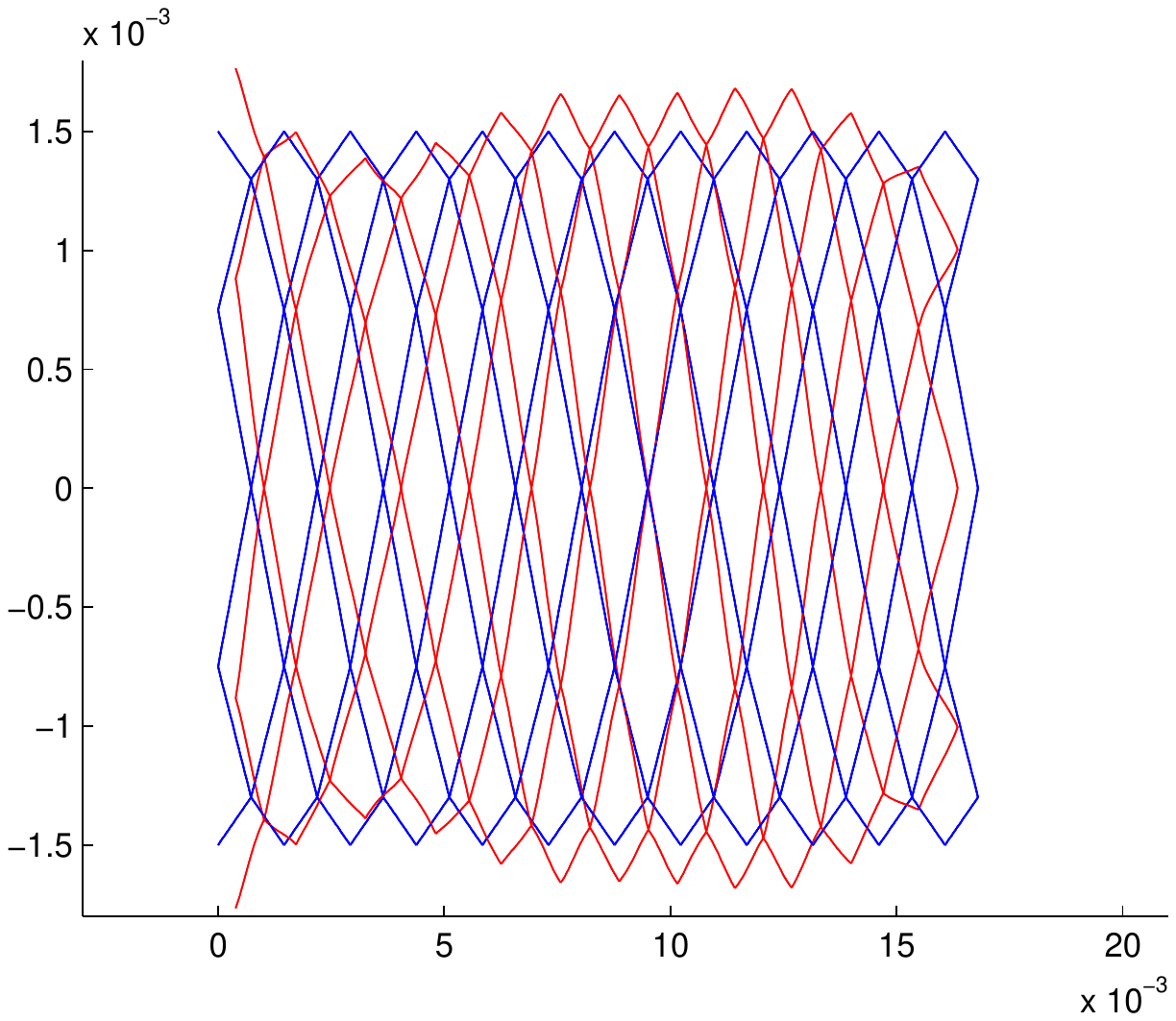}
\caption{Solution of the problem projected to the $(x_1,x_2)$--plane for loads given in (\ref{force}) and at times $t = 1s, 2s, 3s, 4s, 5s$, and $6s$.}\label{time1}
\end{center}
\end{figure}

In Figure~\ref{time1} the solutions of the problem for the force $f$ in \eqref{force} at the time-points $t\in\{1,2,3,4,5,6\}s$ is plotted.

In the sequel we compare the computing  times of two different approaches. In the first approach we use the MATLAB backslash function to solve the system obtained at every time step. In the second, we first perform  the $LDL^T$ decomposition of the matrix $-\Ebb+0.25 \Delta t^2 \Kbb$, since it is the same in all iterations and then in the time integration we use the obtained $LDL^T$ decomposition to solve the system.
%We compare the times needed for the evaluation.

In the Table~\ref{tab3} the computing times for two different calculations with different space discretizations are presented. In the first column the number of splits in the longest edge (strut) in the stent is displayed. The second column displays the associated number of degrees of freedom. The remaining columns present the computing times. The results indicate that the use of the factorization $LDL^T$ obviously pays off.
\begin{table}[h!]\label{tab3}
\begin{center}
\begin{tabular}{c|c|c|c|c}
& & no LDL & \multicolumn{2}{c}{using LDL${}^T$}\\\hline
$\#$ splits & size of matrix & time (s) & time for precomputation (s)& time for iterations (s)\\\hline
$2^{0}$   & 12462          & 476     & 1.34  & 110\\
$2^{1}$   & 25710          & 554     & 1.70  & 203\\
$2^{2}$   & 52206          & 793     & 4.62  & 392\\
$2^{3}$   & 105198         & 1338   &  23.71 & 855\\
\end{tabular}
\end{center}
\caption{Computing times for several space discretizations without and with LDL${}^T$ precomputation.}
\end{table}
Here the time step is equal to $2^{-4}$.

In Table~\ref{tab4} we compare computing times for different time step sizes with the same final time $T$. The total time approximately doubles when lowering $\Delta t$ which is natural, since the number of time steps doubles. However, it is interesting to note that the time for the solution with $LDL^T$ reduces when reducing $\Delta t$.
\begin{table}[h!]\label{tab4}
\begin{center}
\begin{tabular}{c|c|c|c}
& no LDL & \multicolumn{2}{c}{using LDL${}^T$}\\\hline
$\Delta t$ &  time (s) & time for precomputation (s)& time for iterations (s)\\\hline
$2^{-3}$   &  277     & 3.66 & 202 \\
$2^{-4}$   &  554     & 1.77 & 392 \\
$2^{-5}$   &  1120    & 1.52 & 808 \\
$2^{-6}$   &  2332   &  1.39 & 1717\\
\end{tabular}
\end{center}
\caption{Computing times for several values of $\Delta t$ without and with LDL${}^T$ precomputation.}
\end{table}
Here every strut is split in $4$ smaller struts.

We have also computed the errors in the solutions. For the same time step we computed the errors for different number of strut splits. The errors are presented in Table~\ref{errspace}. They are calculated by comparing the solution with the solution for $2^{4}$ strut splits. All errors are presented in the $L^2(0,T;L^2(\cN))$ norm.
\begin{table}[!ht]
\begin{center}
\begin{tabular}{c|c}
$h^{-1}$ & error \\\hline
$2^{0}$   & $0.041110796760794 $ \\
$2^{1}$   & $0.015348501778952 $ \\
$2^{2}$   & $0.003524722458048 $  \\
$2^{3}$   & $0.000224774470258 $\\
\end{tabular}\label{errspace}
\end{center}
\caption{Relative $L^2(0,T;L^2(\cN))$ errors for different space meshes.}
\end{table}
In Table~\ref{errtime} the errors are given for different $\Delta t$ and for fixed space mesh. The errors are computed with respect to $\Delta t = 2^{-7}$ and the same $L^2(0,T;L^2(\cN))$ norm.
\begin{table}[!ht]
\begin{center}
\begin{tabular}{c|c}
$\Delta t$ & error \\\hline
$2^{-3}$   & $0.038370278764979 $ \\
$2^{-4}$   & $0.023014149022735 $ \\
$2^{-5}$   & $0.012051225164378 $ \\
$2^{-6}$   & $0.003154110021893 $\\
\end{tabular}\label{errtime}
\end{center}
\caption{Relative $L^2(0,T;L^2(\cN))$ errors for different time steps.}
\end{table}

Movies of the dynamic behavior of the stent are presented in the video Appendix.
\section{Conclusion}

We have presented a new model description for the numerical simulation of elastic stents, which explicitly displays all constraints. Based on the new formulation an inf-sup inequality for the finite element discretization is shown and, furthermore, a simplified proof of the inf-sup inequality for the space continuous problem is presented. Despite an increased number of degrees of freedom, the new formulation leads to faster simulation time. The  presented techniques  are also  used to simplify the analysis and numerical solution of the evolution problem describing the movement of the stent under external forces. Numerical examples illustrate the theoretical results and show the effectiveness of the new modeling approach.


\begin{thebibliography}{XX}

\bibitem{Bapat}
\textsc{R. B. Bapat}, \textit{Graphs and matrices}, Springer, London; Hindustan Book Agency, New Delhi, 2014.


\bibitem{BeaMXZ18}
\textsc{C. Beattie, V. Mehrmann, H. Xu, and  H. Zwart},
\textit{Linear port-Hamiltonian descriptor systems},
Math. Control Signals Systems, 30 (2018), 17 (27 pages).
%n

%Preprint 06-2017, Institute of Mathematics, TU Berlin, 2017.
%{\it url: http://www.math.tu-berlin.de/preprints/}
%To appear in {\em Math. Control Signals Systems}, 2018.

\bibitem{BBF}
\textsc{D. Boffi, F. Brezzi, and M. Fortin},
\textit{Mixed Finite Element Methods and Applications}, Springer, Heidelberg, 2013.


\bibitem{BrezziFortin}
\textsc{F. Brezzi and M. Fortin},
\textit{Mixed and Hybrid Finite Element Methods},
Springer, New York, 1991.



\bibitem{IMASJ}
\textsc{S. \v{C}ani\'{c} and J. Tamba\v{c}a}, \textit{Cardiovascular stents as PDE nets: 1D vs. 3D}, {IMA J. Appl. Math.}, 77 (2012), pp.~748--770.


\bibitem{CiarletFEM}
\textsc{P.G. Ciarlet}, \textit{The Finite Element Method for Elliptic Problems},  Society for Industrial and Applied Mathematics (SIAM), Philadelphia, 2002.

\bibitem{Evans}
\textsc{L.C. Evans}, \textit{Partial Differential Equations}, American Mathematical Society, Providence, 1998.

\bibitem{GR}
\textsc{V. Girault and P.-A. Raviart}, \textit{Finite Element Methods for Navier-Stokes Equations. Theory and Algorithms},  Springer, Berlin, 1986.


\bibitem{Griso}
\textsc{G. Griso},
\textit{Asymptotic behavior of structures made of curved rods},
Anal. Appl. (Singap.), 6 (2008), pp.~11--22.


\bibitem{RadHAZU}
\textsc{L. Grubi\v{s}i\'{c}, J. Ivekovi\'{c}, J. Tamba\v{c}a, and  B. \v{Z}ugec}, \textit{Mixed formulation of the one-dimensional equilibrium model for elastic stents},  Rad Hrvat. Akad. Znan. Umjet. Mat. Znan., 21(532) (2017), pp.~219--240.

\bibitem{GKMV}
\textsc{L. Grubi\v{s}i\'{c}, V. Kostrykin, K.A. Makarov, and K. Veseli\'{c}}, \textit{Representation theorems for indefinite quadratic forms revisited}, Mathematika, 59 (2013), pp.~169--189.

\bibitem{GT}
\textsc{L. Grubi\v{s}i\'{c} and J. Tamba\v{c}a}, \textit{Direct solution method for the equilibrium problem for
elastic stents}, Numer. Lin. Algebra  Appl., to appear, 2019.

\bibitem{HaiW96}
\textsc{E. Hairer and G. Wanner}, \textit{Solving Ordinary Differential Equations {II}: Stiff and Differential-Algebraic Problems}, 2nd revised edition. Springer Verlag, Heidelberg, 1996.

\bibitem{JT1}
\textsc{M. Jurak and J. Tamba\v{c}a}, \textit{Derivation and justification of
a curved rod model}, Math. Models Methods Appl. Sci., 9 (1999), pp.~991--1014.

\bibitem{JT2}
\textsc{M. Jurak and J. Tamba\v{c}a},
\textit{Linear curved rod model. General curve},
Math. Models Methods Appl. Sci., 11 (2001), 1237--1252.

\bibitem{Kat13}
\textsc{T. Kato}, \textit{Perturbation Theory for Linear Operators}, Springer Science $\&$ Business Media, New York, 2013.

%\bibitem{KT}
%\textsc{M. Kosor and J. Tamba\v{c}a}, \textit{Nonlinear bending-torsion model for curved rods with little regularity}, %Math. Mech. Solids, 22 (2017), pp.~708--717.
%not cited

\bibitem{KunM06}
\textsc{P. Kunkel and V. Mehrmann},
\textit{Differential--Algebraic Equations.
Analysis and Numerical Solution},
European Mathematical Society (EMS), Z{\"u}rich, 2006.

\bibitem{KunMS14}
\textsc{P. Kunkel, V. Mehrmann, and L. Scholz},
\textit{Self-adjoint differential-algebraic equations},
Math. Control Signals Systems, 26 (2014), 47--76.
%n

\bibitem{LamMT13},
\textsc{R. Lamour, R. M{\"a}rz, and C. Tischendorf}, \textit{Differential-Algebraic Equations: A Projector Based Analysis}, Springer Science $\&$ Business Media, Heidelberg, 2013.


\bibitem{MacMMM06b}%101
\textsc{D.S. Mackey, N.~Mackey, C.~Mehl, and V.~Mehrmann},
\textit{Structured polynomial eigenvalue problems:
good vibrations from good linearizations},
SIAM J. Matrix Anal. Appl., 28 (2006), pp.~1029--1051.
%n

\bibitem{MenM03}
\textsc{R. Mennicken and M. M\"{o}ller}, \textit{Non-Self-Adjoint Boundary Eigenvalue Problems}, North-Holland Publishing Co., Amsterdam, 2003.

\bibitem{Tevolution}
\textsc{J. Tamba\v{c}a}, \textit{Justification of the dynamic model of curved rods}, Asympt. Anal.,
31 (2002), 43--68.

\bibitem{SIAMstent}
\textsc{J. Tamba\v{c}a, M. Kosor, S. \v{C}ani\'{c}, and D. Paniagua},
\textit{Mathematical modeling of vascular stents},
SIAM J. Appl. Math., 70 (2010), pp.~1922--1952.


%\bibitem{Zugec1}
%\textsc{J. Tamba\v{c}a and B. \v{Z}ugec},  \textit{A biodegradable elastic stent model}, Math.  Mech. Solids, to %appear, DOI: 10.1177/1081286518773830.

\end{thebibliography}
\end{document}